\input amstex

\documentstyle{amsppt}

\NoRunningHeads \magnification=1200 \TagsOnRight \NoBlackBoxes
\hsize=5.3in \vsize=7.2in \hoffset= 0in \baselineskip=12pt
\topmatter
\title
Convexity of coverings of projective varieties and vanishing
theorems
\endtitle
\rightheadtext{Vanishing theorems}
\author Fedor Bogomolov*
 Bruno De Oliveira**
\endauthor
\thanks
* Partially supported by NSF grant DMS-0100837.
** Partially supported by NSF Postoctoral Research Fellowship
DMS-9902393 and NSF grant DMS-0306487
\endthanks
\affil
* Courant Institute for Mathematical Sciences
**University of Miami
\endaffil
\address
Fedor Bogomolov Courant Institute for Mathematical Sciences,
New York University\\
Bruno De Oliveira University of Miami
\endaddress
\email bogomolo{\@}CIMS.NYU.EDU bdeolive@math.miami.edu
\endemail

\abstract Let $X$ be a projective manifold, $\rho:\tilde X \to X$
its universal covering and $\rho^*: Vect (X) \to Vect(\tilde X)$
the pullback map for isomorphism classes of vector bundles. This
article makes the connection between the properties of the
pullback map $\rho^*$ and the properties of the function theory on
$\tilde X$. Our approach motivates a weakened version of the
Shafarevich conjecture: the universal covering $\tilde X$ of a
projective manifold $X$ is holomorphically convex modulo the
pre-image $\rho^{-1}(Z)$ of a subvariety $Z \subset X$. We prove
this conjecture for projective varieties $X$ whose pullback map
$\rho^*$ identifies a nontrivial extension of a negative vector
bundle $V$ by $\Cal O$ with the trivial extension. We prove the
following pivotal result: if a universal cover of a projective
variety has no nonconstant holomorphic functions then the pullback
map $\rho^*$ is almost an imbedding. Our methods also give a new
proof of $H^1(X,V)=0$ for negative vector bundles $V$ over a
compact complex manifold $X$ whose rank is smaller than the
dimension of $X$.
\endabstract
\endtopmatter

\document

\head 1. Introduction
\endhead

\

This paper deals with two questions about the function theory of
universal covers $\tilde X$ of projective varieties $X$. One
question, appearing in section 2, is on the abundance of
holomorphic functions on $\tilde X$. The main conjecture on the
abundance is the Shafarevich uniformization conjecture (see
below). The other question, appearing in section 3, is on the
simple existence of nonconstant holomorphic functions on $\tilde
X$. It is an open question to know whether the universal cover of
a projective variety has nonconstant holomorphic functions. In
dealing with both problems, we use the same idea. The idea is to
explore the relation between the existence of nonconstant
holomorphic functions on $\tilde X$ and the identification on
$\tilde X$ of the pullback of distinct isomorphism classes of
vector bundles on $X$. This relation gives a new approach to the
production of holomorphic functions on universal cover $\tilde X$
of a projective manifold $X$.  Additionally, using the methods of
section 2 we give a new proof of the vanishing of $H^1(X,V)$ for
negative vector bundles $V$ over a compact complex manifold $X$
whose rank is smaller than the dimension of $X$. The following is
a description of what can be found in this paper.
\par

In section 2, we approach the abundance of holomorphic functions
on the universal covers of projective varieties. In the early 70's
I. Shafarevich proposed the following conjecture on the function
theory of universal covers: The universal cover $\tilde X$ of a
projective variety $X$ is holomorphically convex, i.e every
discrete sequence of points of $\tilde X$ has a holomorphic
function on $\tilde X$ that is unbounded on it. This conjecture
has been proved in some cases (see [Ka95] and [EKPR] for the
strongest results), but the general case has remained unreachable.
For the general case, there is the work of Kollar [Ko93] on the
existence of Shafarevich maps (Campana [Ca94] has dealt with
Kahler case). The existence and properties of the Shafarevich maps
do not give information on the existence of holomorphic functions
on $\tilde X$. But they are an essential tool for dealing with and
understanding the conjecture. There are two main reasons for the
difficulty in proving the Shafarevich conjecture. The first reason
is that the conjecture proposes that noncompact universal covers
$\tilde X$ have many holomorphic functions. But, on the other
hand, there is a lack of methods to construct holomorphic
functions on $\tilde X$. The second reason comes from the main
geometric obstruction to holomorphic convexity. A holomorphic
convex analytic space can not have an infinite chain of compact
subvarieties. The existence of these infinite chains on universal
covers of projective varieties has not been ruled out. In fact,
the first author and L. Katzarkov produced some examples of
algebraic surfaces that possibly contain infinite chains
[BoKa98].\par

The possible existence of infinite chains on universal covers
demands that one rethinks the Shafarevich conjecture. We do
exactly that in section 2, where we use our approach for obtaining
functions on the universal cover to motivate the following
weakened Shafarevich conjecture.

\proclaim {Conjecture 2.3} The universal covering $\tilde X$ of a
projective variety $X$ is holomorphically convex modulo the
pre-image  of a subvariety $Z \subset X$.
\endproclaim

This means that for every infinite discrete sequence $\{x_i\}_{i
\in \Bbb N}$ $x_i \in \tilde X$ such that $\{\rho(x_i)\}$ has no
accumulation points on $Z$, there exists a holomorphic function
$f$ on $\tilde X$ which is unbounded on the sequence. The interest
of this conjecture is that it is still very strong but does not
exclude the existence of infinite chains of compact subvarieties.
The strength of our weakened conjecture is manifested in the fact
that it would still separate universal covers of projective
varieties from universal covers of compact non-kahler manifolds
with many holomorphic functions. The example to have in mind is
the case of the universal cover of an Hopf surface which is $\Bbb
C^2 \setminus \{(0,0)\}$. The complex manifold  $\Bbb C^2
\setminus \{(0,0)\}$ has many holomorphic functions, but it is not
holomorphic convex modulo of the pre-image of any subvariety of
the Hopf surface.\par
\newpage

An explicit motivation for the weakened conjecture can be found
 in theorem A. This theorem  proves the conjecture for
projective varieties $X$ satisfying: $X$ has a negative bundle $V$
such that the pullback map identifies a nontrivial extension of
$\Cal O$ by $V$  the trivial extension.

\proclaim{Theorem A} Let $X$ be a projective variety with a
negative vector bundle $V$ and $\rho:\tilde X \to X$ its universal
covering. If there exists a nontrivial cocycle $s\in H^1(X,V)$
such that $\rho^*s=0$ then $\tilde X$ is holomorphic convex modulo
$\rho^{-1}(Z)$, $Z$ is a subvariety of $X$.
\endproclaim

The nature of the  method  used to produce holomorphic functions
in the proof of the theorem is another motivation for the
conjecture. It will be seen in subsection 2.3 that the method
gives: 1) very strong and precise holomorphic convexity properties
for $\tilde X$; 2) a subvariety $Z$ of $X$ for which the
holomorphic functions on $\tilde X$ created by the method must be
constant over $\rho^{-1}(Z)$. The natural appearance of the
subvariety $Z$ is not one of the method's shortfalls, but rather
one of its strengths, since it is  the existence of $Z$ that
permits the possible existence of infinite chains. If there are
configurations of subvarieties of $X$ whose pre-image contain
infinite chains, they must be contained in $Z$.  To quickly put in
perspective the scope of theorem A, we note that in corollary 2.6
we show that the conditions of theorem A imply that $X$ has a
generically large fundamental group, i.e the general fiber
Shafarevich map is zero dimensional. We note that the varieties
with generically large fundamental group form a natural class of
manifolds to consider when studying the Shafarevich conjecture
[Ko93]. In particular, all the difficulties of the conjecture are
present for this class of manifolds.
\par

In section 3, we deal with the existence of nonconstant
holomorphic functions on the universal cover $\tilde X$ of a
projective variety $X$. The known paths to the production of
holomorphic functions on $\tilde X$ involve the construction of
closed holomorphic 1-forms or exhaustion functions with
plurisubharmonic properties on $\tilde X$. The construction of the
desired closed (1,0)-forms or exhaustion functions on $\tilde X$
involve  the following methods: (a) properties of the fundamental
group $\pi_1(X)$ in combination with Hodge theory and non-abelian
Hodge theory (see [Si88] and [EKPR03] for the most recent results
and references); (b) curvature properties of $X$ (see for example
 [SiYa77] and [GrWu77]), (c) explicit descriptions of $X$ (see for
example [Gu87] and [Na90]). None of these methods are at the
moment sufficiently general to provide a nonconstant holomorphic
function for the universal cover of an arbitrary projective
variety.
\par

Our approach to the existence of holomorphic functions on $\tilde
X$ is different. We connect the existence of nonconstant
holomorphic functions on $\tilde X$ with properties of $\rho^*:
Vect (X) \to Vect(\tilde X)$, the pullback map for vector bundles.
To make our point, we give an extreme example where $\rho^*$
identifies many isomorphism classes. Let $X$ be a projective
manifold such that the pullback map identifies all isomorphism
classes of holomorphic vector bundles on $X$ that are isomorphic
as topological bundles. Then in observation 3.1 we show that
$\tilde X$ must be Stein.
\par

We are interested in obtaining conditions  on the pullback
$\rho^*$ from the absence of nonconstant holomorphic functions on
$\tilde X$. To accomplish this goal, we reexamine the method to
produce functions employed in theorem A. The point is that to able
to obtain a holomorphic function on $\tilde X$ from a cocycle
$s\in H^1(X,V)$ such that $\rho^*s=0$ we do not need $V$ to be
negative. We need, actually, very weak negativity properties on
$V$ as will be illustrated in lemma 3.14 and proposition 3.16.
These results are used to prove the main theorem of this section
describing the pullback map for absolutely stable vector bundles.
A vector bundle is absolutely stable if for any coherent subsheaf
$\Cal F\subset E$ with $\text {rk} \Cal F < \text {rk E}$ a
multiple of the line bundle $(\text {rk E} \text {det} \Cal F
-\text {rk} \Cal F \text {det E})^*$ can be represented by a
nonzero effective divisor. In particular, absolutely stable
bundles are stable with respect to all polarizations of $X$. For
projective surfaces absolutely stable bundles are exactly the
bundles that are stable with respect to all elements in the
closure of the polarization cone.

\proclaim {Theorem B} Let $X$ be a projective manifold whose
universal cover has only constant holomorphic functions. Then:
\par

a) The pullback map $\rho_0^*:Mod_0(X) \to Vect(\tilde X)$  is a
local embedding ($Mod_0(X)$ is the moduli space of absolutely
stable bundles).
\par b) For any absolutely stable bundle $E$ there are
only a finite number of bundles $F$ with $\rho^*E = \rho^* F$.\par
c) Moreover, there is a finite unramified cover $p: X'\to X$
associated with $E$ of degree $d<rk E!$ with universal covering
${\rho'}:\tilde X \to X'$. On $X'$ there is  a collection of
vector bundles $\{E'_i\}_{i=1,...,m}$ on $X'$ with $H^0(\tilde
X,End_0{\rho'}^*E'_i)=0$ such that $\rho^*F\simeq \rho^*E$ if and
only if:
 $$p^*F =E'_1\otimes \Cal O(\tau_1)\oplus ... \oplus E'_1\otimes \Cal
 O(\tau_m)$$
 The bundles $\Cal O(\tau_i)$  are flat bundles associated with finite linear
representations of $\pi_1(X')$ of a fixed rank $k$ with $rk E|k$.
\endproclaim

 The theorem  imposes strong constraints on the pullback map $\rho^*$ if
$\tilde X$ is not to have nonconstant holomorphic functions. In
particular, it says that the pullback map should be almost an
embedding. The authors believe that this imposition on the
pullback map $\rho^*$ should not hold for any projective variety
(see the remarks at the end of section 3). If  the authors are
correct then theorem B would imply the existence of nonconstant
holomorphic functions on the universal cover of a projective
variety.

An interesting application of our method used to prove theorem A
is a new simple proof of the vanishing of the first cohomology
group from negative vector bundles whose rank is smaller than the
dimension of the base.
\medskip

We are grateful to L.Katzarkov, T.Pantev and T. Napier  for their
useful comments. The first author also wants to thank IHES,
University of Bayreuth and University of Miami for partial
support. The second author also thanks the Centro de Analise,
Geometria e Sistemas dinamicos of the IST of Lisbon.

\

\

\

\

\head 2. Convexity properties of universal covers\endhead

\

The first two subsections provide the background and notation for
the subsection 2.3 where the weakened Shafarevich conjecture is
discussed.

\

\subhead {2.1 Holomorphic convexity of universal
covers}\endsubhead

\

A complex manifold $X$ is $holomorphic$ $convex$ if for every
infinite discrete sequence $\{x_i\}_{i \in \Bbb N}$ of points in
$X$ there exists a holomorphic function $f$ on $X$ which is
unbounded on the sequence. Shafarevich proposed the following:

\proclaim {Conjecture}(Shafarevich) The universal cover of a
projective variety is holomorphic convex.
\endproclaim

By the time the conjecture was formulated, it was known that the
holomorphic convexity is a property that is shared by the compact
and noncompact universal covers of Riemann surfaces. It was also
known that if the fundamental group of a Kahler variety $X$ is
abelian then its universal cover $\tilde X$ is holomorphic convex.
But maybe, the most inspiring evidence was the result of Poincare
[Si48] stating that if a bounded domain $D$ in an complex
Euclidean space is the universal cover of a projective manifold,
then $D$ is holomorphic convex.\par

The Shafarevich Conjecture predicts that noncompact universal
covers of projective manifolds have many holomorphic functions.
The holomorphic convexity of the universal cover implies that
there is a proper map of $\tilde X$ into $\Bbb C^n$. In
particular, holomorphic convexity implies that there are enough
holomorphic functions to separate points that are not connected by
a chain of compact analytic subvarieties. This pointwise
holomorphic separability property is the strongest possible for a
complex manifold.
\par

We propose a weakened version of holomorphic convexity, that will
appear in subsection 2.3 to generalize the Shafarevich conjecture.

\proclaim {Definition 2.1} Let $X$ be a complex manifold and
$\rho: \tilde X \to X$ the universal covering of $X$. The
universal cover $X$ is holomorphic convex modulo an analytic
subset $\rho^{-1}(Z)$, $Z\subset X$, if for every infinite
discrete sequence $\{x_i\}_{i \in \Bbb N}$ $x_i \in \tilde X$ such
that $\{\rho(x_i)\}$ has no accumulation points on $Z$, there
exists a holomorphic function $f$ on $\tilde X$ which is unbounded
on the sequence.
\endproclaim

Both holomorphic convexity and holomorphic convexity modulo an
analytic subset of an universal cover imply the abundance of
holomorphic functions on $\tilde X$.

There is also a generalization of the notion of holomorphic
convexity to line bundles other than $\Cal O$. Let $X$ be a
complex manifold and $L$ a line bundle on $X$ with an Hermitian
metric $h$. $X$ is $holomorphic$ $convex$ $with$ $respect$ $to$
$(L,h)$ if for every infinite discrete sequence $\{x_i\}_{i \in
\Bbb N}$ of points in $X$ there exists  a section $s\in H^0(X,L)$
such that the function $|s|_h$ is unbounded on the sequence
(holomorphic convexity is the special case of holomorphic
convexity with respect to the trivial line bundle equipped with
the trivial metric). In subsection 2.3 we use  the following
result from [Na90] concerning holomorphic convexity with respect
to positive line bundles. Let $X$ be a smooth projective variety
and $L$ a positive line bundle on $X$. If $p\gg0$ then  the
universal cover $\tilde X$, $\rho:\tilde X \to X$, is holomorphic
convex with respect to $(\rho^*L^p,h)$, $h$ is any continuous
Hermitian metric on $L^p$.

In the work of Kollar [Ko93] and Campana [Ca94] on the Shafarevich
conjecture, it was shown that every projective (Kahler) manifold
$X$ has a dominant connected rational (meromorphic) map to a
normal variety (analytic space) $Sh(X)$, $sh:X \dasharrow Sh(X)$
such that:

\smallskip

$ \bullet \text { }$ There are countably many closed proper
subvarieties $D_i \subset X$ such that for every irreducible $Z
\subset X$ with $Z \not \subset \bigcup D_i$, one has:
$sh(Z)=\text {point}$ if and only if $\text {im}[\pi_1(\bar Z)\to
\pi_1(X)]$ is finite, $\bar Z$ is the normalization of $Z$.

\smallskip

The map $sh:X \dasharrow Sh(X)$ is called the $Shafarevich$ $map$
and $Sh(X)$ is called the $Shafarevich$ $variety$ $of$ $X$. If the
Shafarevich conjecture would hold, then Shafarevich map would be a
morphism with the property:
\smallskip

$\bullet \text { }$ For every subvariety $Z \subset X$, $
{sh}(Z)=\text {point}$ iff $\text {im}[\pi_1(\bar Z) \to
\pi_1(X)]$ is finite, $\bar Z$ is the normalization of $Z$.
\smallskip

In subsection 2.3 we present the weakened Shafarevich conjecture
and give a motivation for this conjecture for  projective
manifolds $X$ with a $generically$ $large$ $fundamental$ $group$,
i.e $\dim Sh(X)=\dim X$.

\

\

\subhead {2.2  Affine bundles and the negativity of vector
bundles}\endsubhead

\

We recall a construction of affine bundles associated with
extensions of a given vector bundle $V$. We also describe how the
negativity properties of the vector bundle $V$ influence the
function theory of the affine bundle. \par

Let $X$ be a complex manifold and $V$ a vector bundle of rank $r$
on $X$. We will use the common abuse of notation where $V$ also
denotes the sheaf of sections of $V$. An extension of $\Cal O$ by
a vector bundle $V$ is an exact sequence:

$$0\to V \to V_\alpha \to \Cal O\to 0 \tag 2.1$$

There is a 1-1 natural correspondence between cocycles $\alpha\in
H^1(X,V)$ and isomorphism classes of extensions of $\Cal O$ by
$V$.  The extension (2.1) defines an affine bundle, which consists
of the pre-image in $V_\alpha$ of a constant nonzero section of
the trivial line bundle $\Cal O$. This affine bundle is
independent of the choice of the constant nonzero section and is
denoted by $A_\alpha$.  A cocycle $\alpha$ cohomologous to zero
corresponds to the trivial extension $V_\alpha =V\oplus \Cal O$.
We have that the affine bundle $A_\alpha$ is a vector bundle if
and only if (2.1) splits or equivalently if $\alpha$ is
cohomologous to zero. Also recall that an affine bundle is a
vector bundle if and only if the affine bundle has a section. \par

The affine bundle $A_\alpha$ can be described in an alternative
way. Let $E$ be a vector bundle of rank $r$ over $X$, $p:\Bbb
P(E)\to X$ be the $\Bbb P^{r-1}$-bundle over $X$, whose points in
the fiber $\Bbb P(E)_x$ are the hyperplanes in the vector space
$E_x$, $x\in X$. Associated to a surjection $E\to F \to 0$ of
vector bundles there is  an inclusion $\Bbb P(F) \subset \Bbb
P(E)$ of projective bundles. The affine bundle $A_\alpha$ is $\Bbb
P(V_\alpha^*) \setminus \Bbb P(V^*)$, where the inclusion $\Bbb
P(V^*) \subset \Bbb P(V_\alpha^*)$ comes from
 (2.1) dualized.\par

Let $V$ be a vector bundle over a complex manifold $X$. We  recall
Grauert's characterization of negativity  for vector bundles. The
projective bundle $\Bbb P(V)$ has a naturally defined line bundle
$\Cal O_{\Bbb P(V)}(1)$ on it. The line bundle $\Cal O_{\Bbb
P(V)}(1)$ is the quotient $p^*V/F$, where $F$ is the tautological
hyperplane bundle over $\Bbb P(V)$. The vector bundle $V$ is
$Grauert$ $negative$ if the line bundle $\Cal O_{\Bbb P(V^*)}(1)$
is ample.

The negativity properties of a vector  imply complex analytic
properties of the total space  of the vector bundle and of the
associated affine bundles.   Recall that a complex manifold $X$ is
$q$-$convex$ in the sense of Andreotti-Grauert if there is a
$C^{\infty}$ function $\varphi:X \to \Bbb R$ such that outside a
compact subset $K\subset X$:
\smallskip
(i) \quad The subset $\{x\in X: \varphi(x)<c$ with
$c<\sup_X\varphi\}$ is relative compact.\par

(ii) \quad $\varphi|_{X\setminus K}$ is q-convex i.e. the Levi
form $L(\varphi)=\sum_{i,j}\frac{\partial ^2\varphi}{\partial z_i
\partial \bar z_i}dz_i\otimes d\bar z_j$ has at most $q-1$
non-positive eigenvalues at all $x \in X\setminus K$
\smallskip

In subsection 2.3, we will use the following result: if $V$ is a
negative vector bundle  then any affine bundle $A_\alpha$ is
$1$-convex (strongly pseudoconvex). This implies that $A_\alpha$
is holomorphic convex and it has a bimeromorphic morphism onto a
Stein space.

We give a brief proof of the result mentioned above since it plays
an important role in the proof of theorem A and in any of its
generalizations [DeO03]. Let $E \to F$ be a surjection of vector
bundles and $\Bbb P(F)\subset \Bbb P(E)$ be the respective
inclusion of projective bundles. The expression for the normal
bundle $N_{\Bbb P(F)/ \Bbb P(E)}$ is:

$$N_{\Bbb P(F)/ \Bbb P(E)}\simeq p^*(E^*/F^*)\otimes \Cal O_{\Bbb P(F)}(1) \tag 2.2$$

\noindent where $p: \Bbb P(F) \to X$ is the projection.  Consider
the surjection $V_\alpha^*\to V^*$, coming from (2.1) dualized,
and its induced inclusion $\Bbb P(V^*)\subset \Bbb P(V^*_\alpha)$.
The normal bundle $N_{\Bbb P(V^*)/ \Bbb P(V_\alpha^*}\simeq\Cal
O_{\Bbb P(V^*)}(1)$ is positive since $V$ is negative. This
implies that complement $\Bbb P(V_\alpha^*) \setminus \Bbb
P(V^*)=A_s$ is strongly pseudoconvex. As a special case, we have
that if the extension $V_\alpha^*$ is trivial,
$V^*_\alpha=V^*\oplus \Cal O$, then:

$$\Bbb P(V^*\oplus \Cal O) \setminus \Bbb P(V^*)  \simeq V \tag 2.3$$

\noindent Hence the total space $t(V)$ of $V$ is strongly
pseudoconvex.

\

The following is a method to construct many negative bundles of
rank $\geq \dim X$ with nontrivial first cohomology.  Let $L$ be a
very ample line bundle on $X$ which gives an embedding $X \subset
{\Bbb P}^{n}$. There is a surjective map $h: {\Cal O}_{X}^{\oplus
n+1} \to L$ which defines a rank $n$ subbundle $\ker h = F \subset
{\Cal O}_{X}^{\oplus n+1}$. The extension

$$0 \to F\otimes L^{-1}\to \bigoplus^{n+1} L^{-1} \to {\Cal O} \to
0 \tag 2.4$$

\noindent is the pullback  of the Euler exact sequence of $\Bbb
P^n$ to $X$. The vector bundle $F\otimes L^{-1} $ is a negative
bundle, $F\otimes L^{-1} \cong \Omega^{1}_{{\Bbb P}^{n}|X}$, and
$H^1 (X,F\otimes L^{-1})\neq 0$. Namely there is a nontrivial
element $s \in H^1(X,F\otimes L^{-1})\neq 0$ which corresponds to
the above nontrivial extension.

\

The results presented in this paper spring from the existence of
nontrivial  cocycles $\alpha \in H^1(X,V)$ that become trivial
when pulled back to the universal cover. The following standard
result (see [La]) shows that this is only possible for infinite
covers.

\proclaim {Lemma 2.2} Let $f:Y \to X$ a finite morphism between
irreducible normal varieties $X$ and $Y$ and $V$ a vector bundle
over $X$. If $s\in H^1(X,V)$ is nontrivial then $f^*s \in
H^1(Y,f^*V)$ is also nontrivial.
\endproclaim

This result will also be used in theorem A to bound the dimension
of the compact subvarieties in affine bundles associated with
nontrivial cocycles.

\

\

\subhead 2.3 The weakened Shafarevich conjecture\endsubhead

\

The Shafarevich conjecture claims that a noncompact universal
cover of a projective manifold has many holomorphic functions but
it also claims  the  non-existence of infinite chains. As
mentioned in the introduction, the second claim may not hold. In
this subsection we describe an approach to obtain information on
the the algebra of holomorphic functions of universal covers which
has natural place for infinite chains. The main result of this
section, theorem A, motivates the following weakened conjecture:

\proclaim {Conjecture 2.3} The universal covering $\tilde X$ of a
projective manifold $X$ is holomorphically convex modulo the
pre-image of a subvariety $Z \subset X$.
\endproclaim

\noindent As the Shafarevich conjecture, our conjecture also
claims a rich algebra of holomorphic functions but it allows the
existence of infinite chains. The infinite chains of compact
analytic subvarieties would lie in the pre-image of the subvariety
$Z \subset X$ described in the conjecture.
\par

We  describe briefly the methodology our approach. Let $X$ be a
projective manifold, $\rho: \tilde X \to X$ be the universal
covering and $\Cal O_{\tilde X}(\tilde X)$ the algebra of global
holomorphic functions of $\tilde X$. We derive properties of $\Cal
O_{\tilde X}(\tilde X)$ from the existence of nontrivial
extensions of $\Cal O$ by negative vector bundles $V$ which become
trivial once pulled back to $\tilde X$. From such a special
nontrivial extension we construct a map from the universal cover
$\tilde X$ to the associated affine bundle. This map is a local
embedding. From subsection 2.2, it follows that the negativity of
$V$ implies strong analytic geometric properties of the algebra of
global holomorphic functions of the affine bundle. We use these
analytic geometric properties and the local embedding of $\tilde
X$ in the associated affine bundle to obtain information on $\Cal
O_{\tilde X}(\tilde X)$.\par

Our approach motivates the conjecture in two levels. First, it
gives in theorem A an explicit confirmation of the conjecture for
projective manifolds $X$ having  a nontrivial extension of $\Cal
O$ by a negative vector bundle $V$ whose pullback to $\tilde X$ is
trivial. Second, it is the nature of the approach to give very
strong holomorphic convexity properties for $\tilde X$ but also to
give a subvariety $Z$ of $X$ for which all the holomorphic
functions on $\tilde X$ created by the method must be constant
over $\rho^{-1}(Z)$. The subvariety $Z$ is the projection into $X$
of the maximal compact analytic subset of the affine bundle
associated with the extension (see the proof of theorem A).
\par

\proclaim{Theorem A} Let $X$ be a projective manifold with a
negative vector bundle $V$ and $\rho:\tilde X \to X$ its universal
covering. If there exists a nontrivial cocycle $s\in H^1(X,V)$
such that $\rho^*s=0$ then $\tilde X$ is holomorphic convex modulo
$\rho^{-1}(Z)$, $Z$ is a subvariety of $X$.
\endproclaim

\demo {Proof} First, we will identify the subvariety $Z\subset X$
described in the theorem. The nontrivial cocycle $s\in H^1(X,V)$
has an associated strongly pseudoconvex affine bundle $A_s=\Bbb
P(V_s^*) \setminus \Bbb P(V^*)$ originating from the nonsplit
exact sequence:

$$ 0\to \Cal O_X\to V^*_s \to V^* \to 0 \tag 2.5$$

\noindent The  strongly pseudoconvex manifold $A_s$ (hence
holomorphic convex) has a proper holomorphic map onto a Stein
space, $r:A_s \to St(A_s)$ (the Remmert reduction). Moreover,
$A_s$ has a subset $M$ called the maximal compact analytic subset
of $A_s$ such that the map $r|_{A_s\setminus M}: A_s\setminus M
\to r(A_s\setminus M)$, is a biholomorphism. The subvariety
$Z\subset X$ is $Z=p(M)$.\par

\

To obtain holomorphic functions on $\tilde X$ we will construct a
holomorphic map $g: \tilde X \to A_s$ such that $g(\tilde X)\not
\subset M$ and pullback the holomorphic functions of $A_s$ to
$\tilde X$. The pullback of the exact sequence (2.5) to the
universal covering $\tilde X$ splits into:

$$0\to \Cal O_{\tilde X}\to \Cal O_{\tilde X} \oplus \rho^{*} V^* \to \rho^{*} V^* \to 0 $$

\noindent  since it is associated with the trivial cocycle
$\rho^{*}s \in H^1(\tilde X, \rho^*V)$. As observed in (2.3)
$A_{\rho^*s}\equiv \Bbb P(\Cal O_{\tilde X} \oplus \rho^{*} V^*)
\setminus \Bbb P(\rho^*V^*)\simeq \rho^*V$, hence $\tilde X$
embeddeds in $A_{\rho^*s}$ as the zero section of $\rho^*V^*$. The
affine bundle $A_{\rho^*s}$ is the fiber product
$A_{\rho^*_s}=\tilde X \times_X A_s$, denote the projection to the
second factor by $\rho ':A_{\rho^*s} \to A_s$ and   the embedding
of $\tilde X$ in $A_{\rho^*s}$ as the zero section of the vector
bundle $\rho^*V^*$ by $s:\tilde X \to A_{\rho^*s}$. The
holomorphic map $g:\tilde X \to A_s$ will be the composition
$g=\rho' \circ s:\tilde X \to A_s$. The map $g$ is a local
biholomorphism between $\tilde X$ and $g(\tilde X)$ hence the
condition  $g(\tilde X)\not \subset M$ will hold if $\dim M<dim
X$. \par

The maximal compact analytic subset of $A_s$ is of the form
$M=\cup^{k}_{i=1}M_i$, where the $M_i$ are the compact irreducible
positive dimensional subvarieties of $A_s$. The following
proposition shows that  $\dim M_i<dim X$.

\

\proclaim {Proposition 2.5} Let $X$ be a  projective manifold with
a vector bundle $V$ and $V_s$ be the extension associated with a
nontrivial cocycle $s\in H^1(X,V)$. Then any compact subvariety
$M$ of the affine bundle $A_s=\Bbb P(V^*_s) \setminus \Bbb P(V^*)$
satisfies $ {\dim} M< {\dim} X$.
\endproclaim

\demo{Proof} It is clear that if $M \subset A_s$ is a compact
subvariety then $\dim M \le \dim X$ (The intersection of $M$ with
any fiber of the projection map $p: A_s \to X$ will be at most
0-dimensional).  We will show  that $p^*s|_M \in H^1(M,p^*V|_M)$
must be trivial and  that if ${\dim} M={\dim} X$  then $p^*s|_M$
must be nontrivial. These two results prove the desired  strict
inequality $ {\dim} M< {\dim} X$.

The triviality of $p^*s|_M$ follows from the triviality of $p^*s
\in  H^1(A_s,p^*V)$. The equality $p^*s=0$ holds if and only if
the pullback of the exact sequence (3.1) to $A_s$ splits. The
affine bundle $A_s$ is $\Bbb P(V_s^*)\setminus \Bbb P(V^*)$ hence
there is a canonical association between the  points $y\in A_s$
with $p(x)=y$ and the hyperplanes of $(V^*_s)_x$ surjecting to
$(V^*)_x$. From this association one obtains  the canonical
subbundle $W\subset p^*V^*_s$ on $A_s$. The vector bundle $W$ is
such that that the restriction of the surjection $q:V^*_s \to V^*$
to $W$ is still a surjection and hence an isomorphism. The
splitting of $0 \to \Cal O \to p*V^*_s \to p^*V \to 0$ is obtained
by inverting $q:W \to V^*$. If ${\dim} M= {\dim} X$ then $M$ has
an irreducible  component $M'$ such that $p^*s|_{M'}=0$ and the
map $p|_M': M' \to X$ is finite map.  Let $n:\hat M' \to M'$ be
the normalization map, then $n\circ p|_{M'}:\hat M' \to X$ is
finite map between normal varieties hence by the lemma 2.2
$p|_{M'}^*s\neq 0$ which is a contradiction.\hfill \ \qed
\enddemo

Let $\{x_i\}_{i\in \Bbb N}$ be a sequence of points in $\tilde
X\setminus \rho^{-1}(Z)$ such that $\{\rho(x_i)\}_{i \in \Bbb N} $
has no accumulation points on $Z$. The sequence $\{x_i\}_{i\in
\Bbb N}$ has a subsequence $\{y_i\}_{i\in \Bbb N}$ satisfying
$\{\rho(y_i)\}_{i\in \Bbb N}$ converges to $a\in X\setminus Z$.
Consider the sequence $\{g(y_i)\}_{i \in \Bbb N} $ of points in
$A_s$, we have two cases: 1) $\{g(y_i)\}_{i \in \Bbb N} $ is a
discrete sequence of points of $A_s$; 2) $\{g(y_i)\}_{i \in \Bbb
N}$ has a subsequence converging to a point in $a' \in p^{-1}(a)$.
If case 1) holds, since $A_s$ is holomorphic convex it follows
that there is a function $f'\in \Cal O_{A_s}(A_s)$ that is
unbounded on $\{g(y_i)\}_{i \in \Bbb N}$. Hence $f'\circ g$ is the
desired unbounded function on $\{x_i\}_{i\in \Bbb N}$. We procced
to deal with case 2).

Let $L$ be a positive line bundle on $X$, Napier's result [Na90]
states that for $p\gg 0$ $\exists$ $s\in H^0(\tilde X,\rho^*L^p)$
such that $\{|s(y_i)|_{\rho^*h^p}\}_{i\in \Bbb N}$ is unbounded
($h$ is an $C^\infty$ Hermitean metric on $L$). Let $s'\in
H^0(X,L^p)$ be such that $a \notin D=(s')_0$. The meromorphic
function $h=\frac {s}{\rho^*s'}$ is holomorphic outside
$\rho^{-1}(D)$ and unbounded on $\{y_i\}_{i\in \Bbb N}$. Assume
the existence of  a $q \in \Cal O_{\tilde X}(\tilde X)$ with
$\inf\{|q(z)|: z \in g^{-1}(a')\}\neq 0$ and vanishing on
$\rho^{-1}(D)$. Then for $l$ sufficiently large $f=hq^l$ would be
the desired holomorphic function.

An holomorphic function $q \in \Cal O_{\tilde X}(\tilde X)$, as
desired above, can be obtained by pulling back, using $g$, a
holomorphic function $q' \in \Cal O_{A_s}(A_s)$ that satisfies
$q'(p^{-1}(D))=0$ and $q'(a')=1$. The existence of such $q'$
follows from $r:A_s \to St(A_s)$ being a proper map, $r(a')\cap
r(p^{-1}(D \cup Z))=\emptyset$ and $St(A_s)$ being Stein. \hfill \
\qed
\enddemo

\

It is important to complement theorem A with an example that shows
that the hypothesis of the theorem do not imply that the universal
cover $\tilde X$ is Stein. The Steiness of $\tilde X$ holds if the
affine bundle $A_s$ is Stein (corollary 2.7), but otherwise
$\tilde X$ need not be Stein.

\

\noindent Example: We give an example of a projective variety and
a vector bundle satisfying the hypothesis of theorem A but whose
universal cover is not Stein. Let $X$ be a nonsingular projective
variety whose universal cover $\tilde X$ is Stein. Let $\sigma:Y
\to X$ be the blow up of $X$ at a point $p \in X$
$E=\sigma^{-1}(p)=\Bbb P^{n-1}$, $n=\dim X$. If we pullback the
exact sequence (2.5) to $Y$ and tensor it with $\Cal O(E)$ we
obtain:

$$0 \to \sigma^*(F\otimes L^{-1})\otimes \Cal O(E) \to \bigoplus^{n+1}
\sigma^*L^{-1}\otimes \Cal O(E) \to {\Cal O(E)} \to 0 \tag 3.2$$

The vector bundle $\sigma^*(F\otimes L^{-1})\otimes \Cal O(E)$ on
$Y$ is negative since $\sigma^*L^{-1}\otimes \Cal O(E)$ is a
negative line bundle on $Y$. Tensoring a negative vector bundle
with a globally generated vector bundle gives a negative vector
bundle [Ha66]. The pair $Y$ and $\sigma^*(F\otimes L^{-1})\otimes
\Cal O(E)$ satisfies the conditions of theorem A but $\tilde Y$ is
clearly not Stein (it contains $\pi_1(X)$ copies of $\Bbb
P^{n-1}$).\par

\

The existence of a vector bundle $V$ on $X$ satisfying the
hypothesis of theorem A does impose conditions on $X$ and $\tilde
X$. These are described in the following corollary:

\proclaim {Corollary 2.6} Let $X$ be a projective manifold with a
negative vector bundle $V$ and $\rho:\tilde X \to X$ its universal
covering. If it exists a nontrivial cocycle $s\in H^1(X,V)$ such
that $\rho^*s=0$ then:

a) $X$ has a generically large fundamental group.

b) The holomorphic functions on $\tilde X$  separate points on
$\tilde X \setminus \rho^{-1}(Z)$, $Z$ a subvariety of $X$. There
is a finite collection of holomorphic functions such that
$(f_1,...,f_l): \tilde X \to \Bbb C^l$ is a local embedding of
$\tilde X$ at every point in $\tilde X \setminus \rho^{-1}(Z)$.
\endproclaim

\demo {Proof} In the proof of theorem A   a holomorphic map
$g:\tilde X \to A_s$ giving a local embedding of $\tilde X$ into $
A_s$ was constructed. The affine bundle $A_s$ is strongly
pseudoconvex, hence it has the Remmert reduction map $r:A_s \to
St(A_s)$ which is a bimeromorphic morphism. The maximal compact
subset $M$ of $A_s$ consists of the union of the positive
dimensional fibers of the Remmert reduction. Let $Z$ be the
projection of $M\subset A_s$ into $X$. Let $(f'_1,...,f'_l):
St(A_s) \to \Bbb C^l$ be an embedding of the Stein space $St(A_s)$
into an complex Euclidean space. The collection of the functions
$f_i=f'_i\circ r\circ g$ give the local embedding in described b)
since $g(\tilde X)\subsetneq M$ .
\par

To establish b) we still need to show that if $a$ and $b$ are two
different points in $\tilde X \setminus Z$  then there is an $f\in
\Cal O_{\tilde X}(\tilde X)$ such that $f(a)\neq f(b)$. The $\bar
\partial$-method gives  that for any positive line bundle $L$ on $
X$  the line bundle $\rho^*L^p$  on $\tilde X$ for $p\gg 0$ has a
section $s$ with $s(a)=0$ and $s(b)\neq 0$  [Na90]. Let $s'\in
H^0(X,L^p)$ be such that $\rho(b) \notin D=(s')_0$. The
meromorphic function $h=\frac {s}{\rho^*s'}$ is holomorphic
outside $\rho^{-1}(D)$ and $h(b)\neq 0$.  As shown in the proof of
theorem 3.1 there is  a $q \in \Cal O(\tilde X)$ vanishing on
$\rho^{-1}(D)$ and not vanishing at $\rho^{-1}(b)$. Then for $l$
sufficiently large $f=hq^l$ would be the desired holomorphic
function.\par

The projective variety has a generically large fundamental group
since no compact subvariety of $\tilde X$ passes through any point
in $\tilde X \setminus \rho^{-1}(Z)$.

 \hfill \ \qed

\enddemo

\

This corollary defines the natural setting for theorem A. As
mentioned in the introduction,  projective varieties with
generically large fundamental group form a natural class to test
and explore the Shafarevich conjecture. If one wants to obtain a
result similar to theorem A for a general $X$ one needs to
consider semi-negative vector bundles on $X$. This will be done in
the paper [DeO03]. The extensions of semi-negative vector bundles
$V$ associated with a nontrivial $s\in H^1(X,V)$ satisfying
$\rho^*s=0$ exist if the Shafarevich conjecture holds. They can be
obtained as the pullback to $X$ by the Shafarevich map of the
nontrivial extension given by the Euler sequence associated with
an embedding of the Shafarevich variety (the target of the
Shafarevich map).

\

Let us consider the special case where the affine bundle $A_s$
over $X$ is a Stein manifold. The condition that $A_{s}$ is a
Stein manifold can be easily fulfilled in the following examples.
Over $\Bbb P^n$ we have that the  the affine bundle $A_\omega$
associated with the extension $0 \to \Omega^{1}_{{\Bbb P}^{n}} \to
\bigoplus^{n+1} \Cal O(-1) \to {\Cal O} \to 0$ is isomorphic to
the affine variety:

$$F={\Bbb P}^{n}\times {\Bbb
P}^{n \vee}\setminus \{(x,h) \in {\Bbb P}^{n}\times {\Bbb P}^{n
\vee}| x \in h \}
$$

Let $X$ be a projective variety embedded in $\Bbb P^n$ and
$A_\omega|_X$ the affine bundle associated with  pullback to $X$
of the above extension. The affine bundle  $A_\omega|_X$ is a
Stein manifold since it is a closed subvariety of $F$. The
following is a corollary of theorem A for the case where $A_s$ is
Stein.

\proclaim{Corollary 2.7} Let $X$ be a projective manifold, $V$ a
vector bundle and $s\in H^1(X,V)$. Assume furthermore that
$A_{s}=\Bbb P(V^*_s)\setminus \Bbb P(V^*)$ is a Stein variety. Let
$f : Y \to X$ be any infinite unramified covering s.t. $f^{*}s =
0$. Then $Y$ is Stein.
\endproclaim
\demo{Proof} Since any non-ramified covering of a Stein space is
Stein \cite{5} the assumption that $A_{s}$ is affine yields that
$A_{s}\times_{X} Y$ is Stein. On the other hand in the proof of
Theorem A we saw that $Y \subset A_{s}\times_{X} Y$ is a closed
analytic subset and so $Y$ is Stein. \hfill \ \qed
\enddemo

Corollary 2.7 suggests that the result of Theorem A may also be
applicable to orbicoverings of $X$. Let us first describe
precisely the notion of orbicovering in the case of a complex
variety. Let $X$ be a complex variety and $S \subset X$ be a
proper analytic subset. Consider for any point $q\in S$ the local
fundamental group $\pi_q = \pi_1(U(q)\setminus S)$ where $U(q)$ is
a small ball in $X$ centered at $q$. Let $L \subset
\pi_1(X\setminus S)$ be a subgroup with the property that $L\cap
\pi_q$ is of finite index in $\pi_{q}$ for all $q \in S$. Then the
nonramified covering of $ X\setminus S$ corresponding to $L$ can
be naturally completed into a normal complex variety $Y_{L}$ with
a locally finite and locally compact surjective map $f_L : Y_{L}
\to X$. The map $f_L : Y_{L} \to X$ is called an {\it
orbicovering} of $X$ with a ramification set $S$.  The following
holds:

\proclaim{Corollary 2.8}Let $X$ be a projective manifold, $V$ a
vector bundle and $s\in H^1(X,V)$. Assume furthermore that
$A_{s}=\Bbb P(V^*_s)\setminus \Bbb P(V^*)$ is an affine variety.
Let $f : Y \to X$ be any orbicovering s.t. $f^{*}s = 0$. Then $Y$
is Stein.
\endproclaim
\demo{Proof} Since every orbicovering of a Stein space is also
Stein (see Theorem~4.6 of \cite{5}) the proof is exactly the same
as the proof of Corollary 3.5. \hfill \ \qed
\enddemo

\

\

  \head 3. The pullback map for vector bundles and holomorphic functions
on universal covers\endhead

\

\

The positive results on Shafarevich conjecture generally  involve
the existence of nonisomorphic vector bundles on $X$ which become
isomorphic after the pullback to $\widetilde{X}$. For example, the
theorem of L.Katzarkov \cite{4} establishes the holomorphic
convexity of $\widetilde{X}$ for a projective surface $X$ under
the assumption of the existence of an almost faithful linear
representation of $\pi_{1}(X)$. In this case all the bundles on
$X$ corresponding to the representations of the same rank of the
fundamental group  are becoming equal on $\widetilde{X}$. Let $X$
be a projective manifold, $\rho:\tilde X \to X$ its universal
covering and $\rho^*: Vect (X) \to Vect(\tilde X)$ the pullback
map for isomorphism classes of holomorphic vector bundles. This
section investigates the relation between the properties of the
pullback map $\rho^*$ and the existence of holomorphic functions
on $\tilde X$.

We start by considering the case of projective manifolds whose
pullback map $\rho^*$ identifies the isomorphism classes that are
isomorphic as topological vector bundles.  If the pullback satisfy
this property, then there are plenty of distinct vector bundles on
$X$ whose pullbacks are identified. In particular, any two bundles
which can be connected by an analytic deformation are bound to be
identified on $\tilde X$. This very rich collection of bundles
that are identified via the pullback map imply the following
result on the algebra of global holomorphic functions on $\tilde
X$ holds:

\proclaim {Observation 3.1} Let  $X$ be a projective manifold
whose pullback map $\rho^*$ identifies  isomorphism classes of
holomorphic vector bundles that are  in the same topological
isomorphism class. Then the universal cover $\tilde X$
 is Stein.
\endproclaim

\demo {Proof} Let $X$ be a subvariety of $\Bbb P^n$. Let $\alpha
\in H^1(X,\Omega^1_{\Bbb P^n}|_X)$ be the cocycle associated with
the extension of $\Omega^1_{\Bbb P^n}|_X$ coming from the Euler
exact sequence, $(\Omega^1_{\Bbb P^n}|_X)_\alpha$. The pullback
$\rho^*\alpha=0$, since $\rho^* (\Omega^1_{\Bbb P^n}|_X)_\alpha$
is isomorphic to the pullback $\rho^*(\Omega^1_{\Bbb P^n}|_X
\oplus \Cal O)$, topologically they are the same bundle. The
result then follows from corollary 2.7 and the paragraph preceding
it. \hfill \ \qed
\enddemo

The condition that $\rho^*$ identifies isomorphism classes of
holomorphic vector bundles that are  in the same topological
isomorphism class could be replaced by the following apparently
weaker condition:  any extensions of a vector bundle $V$ by
another vector bundle $V'$ are identified under $\rho^*$.

\

This section is mainly concerned with the implications of the
absence of nonconstant holomorphic functions on $\tilde X$ on the
 pullback map $\rho^*$. The condition that $\tilde X$ has no nonconstant
holomorphic functions lies on the opposite side of the conclusion
of observation 3.1, stating that $\tilde X$ is Stein. We will show
that this condition on $\tilde X$ has implications that are quite
opposite to the assumption of the observation 3.1. More precisely,
the absence of nonconstant holomorphic function on $\tilde X$
implies that the pullback map $\rho^*$  is almost an imbedding.
This conclusion lies in strict contrast with the assumption of
3.1, which implies that the pullback map identifies many bundles.
The authors believe that this result should not hold for
projective varieties with  infinite $\pi_1(X)$; see the remarks at
the end of this section. If the authors are correct, our approach
will show that there are always nonconstant holomorphic functions
on noncompact universal covers $\tilde X$ of projective varieties.

Let $\rho^*:Vect(X) \to Vect (\tilde X)$ be the pullback map
sending the set of isomorphism classes of vector bundles on $X$
into the set of isomorphism classes of vector bundles on $\tilde
X$. The flat vector bundle on $X$ obtained from a linear
representation $\tau$ of the fundamental group of $X$ or its sheaf
of sections is denoted by $\Cal O(\tau)$. By construction
$\rho^*\Cal O(\tau)$ is the trivial bundle on the universal
covering $\tilde X$ with the rank of $\tau$. It is clear that two
bundles $F\otimes \Cal O(\tau)$ and $F\otimes \Cal O(\tau')$
become isomorphic on $\tilde X$ for any bundle $F$ if the rank of
the representations $\tau$ and $\tau'$ is the same. The main
result of this section states that if $\tilde X$ has no
holomorphic functions and two bundles $E$ and $E'$ on $X$ have
isomorphic pullback on $\tilde X$ then $E=F\otimes \Cal O(\tau)$
and $E'=F \otimes \Cal O(\tau')$ for some bundle $F$.

In order to better understand the map $\rho^*$ and, in particular,
to study its local properties, one should put a structure of an
analytic scheme into the sets $\text {Vect}(X)$ and $\text
{Vect}(\tilde X)$. The  scheme structure for $\text {Vect}(X)$ for
$X$ projective is well understood. Below, we recall the key facts
that are  relevant to our goals. The analytic scheme structure
theory for $\text {Vect} (\tilde X)$ is less understood. We would
like to note that it is not our interest to develop such a theory
in this paper. The only facts that we use from $\text {Vect}
(\tilde X)$ are that distinct points correspond to non-isomorphic
bundles and that the formal tangent space at a vector bundle $E$
on $\tilde X$ exists and it is equal to $H^1(\tilde X, End E)$.

\

\subhead 3.1 Stability Background \endsubhead

\

To obtain a good parameterizing scheme for vector bundles on a
projective variety $X$ we have to consider some stability
conditions on the bundles, see below. There is an algebraic
parameterization for $H$-stable bundles with given topological
invariants. This parameterization space has all the basic
properties of a coarse moduli space (see for example [HuLe97],
[Ma77]).

Let $E$ be a vector bundle on a projective variety $X$ of
dimension $n$ and $H$ be an arbitrary element in the closure of
the polarization cone $P\subset H^{1,1}(X,\Bbb R)$. $E$ is said to
be $H-semistable$ if   the inequality $(\text {rk E} det \Cal F
-\text {rk} \Cal F det E). H^{n-1}\leq 0$ holds for all coherent
subsheaves $\Cal F\subset E$. Moreover, if for all coherent
subsheaves $\Cal F\subset E$ of lower rank  $(\text {rk E} det
\Cal F -\text {rk }\Cal F det E). H^{n-1} < 0$ holds then $E$ is
said to be $H-stable$. The vector bundle  $E$ is $H-unstable$  if
it has an $H$-destabilizing subsheaf $\Cal F$, i.e  there is a
coherent subsheaf $\Cal F$ with $0<\text {rk} \Cal F<\text {rk E}$
such that $(\text {rk E} det \Cal F -\text {rk} \Cal F det E).
H^{n-1}
> 0$ holds. The number $\mu_H(\Cal F)=(det F/\text {rk} \Cal F). H^{n-1}$ is called
the $H-slope$ of $F$. $H$-stability of $E$ is equivalent to the
fact that any coherent subsheaf of $E$ with smaller rank has a
 smaller $H$-slope than $E$. The notion of $H$-stability for $H$
 is the same as for $aH$, $a\in
\Bbb R$ and $a>0$.

Denote by $P(X)$ the polarization cone of $X$ in the real space
$H^{1,1}(X,\Bbb R)$ and denote its closure by $\bar P(X)$. Since
the base of the closure $\Bbb P(\bar P(X))$ of the cone $\bar
P(X)$ in real projective space $\Bbb P (H^{1,1}(X,\Bbb R))$ is
compact the notion of a stable bundle with respect to the closure
of the polarization cone is well defined. In particular we have
the following result (see [Bo78], [Bo94] and [HuLe97]).

\

\proclaim {Definition 3.2}  The cone of effective
 divisors on $X$, denoted by $K_{eff}$,
 is the cone in the group $\text{Pic (X)}\otimes \Bbb R$
 generated using only non-negative real coefficients by the representatives of effective
 divisors.
 The cone of effective divisors without the zero, $K_{eff} \setminus
 0$, will be denoted by $K_{eff^+}$. Similarly the cone of
anti-effective divisors is denoted by $- K_{eff}$ and
$-K_{eff^+}$.
\endproclaim

\proclaim{Lemma 3.3} Let $E$ be a vector bundle over projective
surface $X$ which is stable with respect to all elements $H$ in
$\bar P(X)$. Then any coherent proper subsheaf $\Cal F\subset E$
satisfies $rk F det E - rk E det \Cal F \in K_{eff^+}$.
\endproclaim

\demo{Proof} Thanks to the stability assumption we have that $(rk
\Cal F det E - rk E det \Cal F).H > 0$ for any $H$ in the closure
of the polarization cone. Since the cone $\Bbb P(\bar P(X))$ is
compact, it is also true in the neighborhood of the cone. Thus,
using Kleiman duality for surfaces,  we obtain that a positive
multiple of $(\text {rk} \Cal F \text {det E} - \text {rk E} \text
{det} \Cal F)$ is effective and nonzero.\enddemo

For a general projective variety $X$ the stability property is
captured on surfaces which are complete intersections in the
initial variety. Thus we make the following definition (see also
[Bo78])

\proclaim{Definition 3.4} A vector bundle $E$ is absolutely stable
if for any coherent subsheaf $\Cal F\subset E$ with $\text {rk}
\Cal F < \text {rk E}$ the following holds: $ \text {rk E} \text
{det} \Cal F -\text {rk} \Cal F \text {det E}$ belongs to the cone
$-K_{eff^+}$. In particular an absolutely stable bundle is stable
with respect to all polarizations.
\endproclaim

\

 The condition of absolute stability is the right
condition for the formulation of our results later in this
section. We will need some properties of the bundle $\text {End}
E$ for an $H$-stable bundle $E$. We will also need results on the
theory of stable bundles for smooth projective curves and on how
stability behaves under restriction maps. We start with the basic
lemma (for a proof see, for example, chapter 1 of [HuLe97]):

\proclaim{Lemma 3.5} Let $E$ be a vector bundle on $X$ which is
stable with respect to some $H\in P(X)$. Then $H^0(X,End_0E) = 0$.

\endproclaim

If $C$ is a smooth curve and $E$ is a stable vector bundle over
$C$ then by a classical result of Narasimhan-Seshadri $End _0E$ is
obtained from a unitary representation $\tau$ of the fundamental
group $\pi_1(C)$ in $PSU(n),n = rk E$. The elements of $PSU(n)$
act on the  matrices in $\text {End} \Bbb C^n$ by conjugation.
Since the bundle $E$ is stable the representation $\tau$ is
irreducible.

\proclaim{Lemma 3.6} If $E$ is a stable vector bundle over a
smooth projective curve $C$ then the bundle $End E$ is a direct
sum of stable vector bundles of degree $0$.
\endproclaim

\demo{Proof} This fact is well known and the decomposition into a
direct sum of stable bundles corresponds to the decomposition of
the unitary representation of $\pi_1(X)$ in $PSU(n)\subset SU(n^2
-1)$ under the above imbedding.\hfill \ \qed
\enddemo

Let $E$ be a vector bundle over a projective variety $X$. Let $C_E
\subset Pic X\otimes \Bbb R$ be the cone generated by the classes
$det \Cal L$ for all coherent subsheaves $\Cal L\subset \bigcup
_{n=0}^\infty(End E)^{\otimes n}$ with $\text {rk } \Cal L=1$. Let
us remember the following result from [Bo94] which follows from
invariant theory:

\proclaim {Lemma 3.7} The cone $C_E$ is also generated  by the
elements of the form: \newline \noindent $det\Cal F - (rk \Cal
F_i/rk E)det E $, where $\Cal F\subset E$ is a proper coherent
subsheaf of $E$ and some elements in $-K_{eff}$. For any bundle
$E$ of rank $k$ there is a natural reductive structure group
$G_E\subset GL(k)$ of $E$ such that $C_E$ is generated by the line
subbundles $L$ corresponding to the characters of parabolic
subgroups in $G_E$.
\endproclaim

The group $G_E$ is defined modulo scalars by the set of subbundles
$L \in (End E)^{\otimes n}$ for all $n$ with $c_1(L) = 0$. If $G_E
= GL(k),SL(k)$ then the line subbundles $L$ are exactly the line
bundles $det\Cal F - (rk \Cal F_i/rk E)det E $. However, if the
group $G_E$ is smaller than above then the line bundles generating
$C_E$ correspond to determinants of special subsheaves of $E$.

\proclaim{Corollary 3.8} Let $E$ be an absolutely stable vector
bundle on a smooth projective variety $X$ and $\Cal A\subset End
E$ be a coherent subsheaf. Then $det \Cal  A \in -K_{eff}$. If $E$
is absolutely stable and $E= E'\otimes F$  then both $E',F$ are
absolutely stable since the corresponding parabolic group $G_E$ is
contained in the  group product $G_{E'}\times G_F$ and the cone
$C_E$ is a sum $C_{E'} + C_F$.

\endproclaim

\demo {Proof} Since $det \Cal A \subset (End E)^{\otimes \text
{rk} End E}$ , $det \Cal A \in C_E$ and the previous lemma implies
that $det \Cal A= \Sigma a_i(det \Cal F_i - (rk \Cal F_i/rk E)det
E)$ where $a_i\geq 0$ and $\Cal F_i$ coherent subsheaves of $E$.
The conclusion follows from the condition of absolute stability,
i.e all elements $det \Cal F_i - det E (rk \Cal F_i/rk E)$ belong
to $-K_{eff}$.
 The parabolic subgroups in $G_{E'}\times G_F$ are products
 of the parabolic subgroups in $G_{E'}, G_F$ which implies
 the result.

\hfill \ \qed
\enddemo

Many properties of $H$-stable bundles on arbitrary projective
varieties can be derived from their restrictions on smooth curves.
As is manifested in the following two results.

\proclaim{Lemma 3.9} Let $X$ be a  projective variety, $H$ a
polarization of $X$, $E$ an $H$-stable vector bundle and $C$ a
generic curve in $kH^{n-1}$ for $k\gg 0$. Then:

1) The restriction of $E$ to $C$ is stable.

2) Any saturated coherent subsheaf $\Cal F\subset EndE|_C$ with
$\mu_H(\Cal F)=0$ is a direct summand of $End E|_C$.

3) The set of saturated subsheaves $\Cal F$ of $End E$ with
$\mu_H(\Cal F) = 0$ coincide via the restriction map with the
similar set for $End E|C$ on $C$.

4) The bundle $End E$ is $H$-semistable and it is a direct sum of
$H$-stable bundles $F_i$ with $\mu_H(F_i) = 0$.

\endproclaim

\demo{Proof}   1), 2) follows from general results, see for
example [Bo78], [Bo94] and [HuLe97]. 3) is a consequence of the
following result: the algebras $\bigoplus^{\infty}_{i=1} H^0(X,
S^i End_0 E)$ and $\bigoplus^{\infty}_{i=1} H^0(C, S^i End_0
E|_C)$ are isomorphic up to a level $l$ depending on $C$
(depending on $k$). This isomorphism follows from the vanishing of
the cohomology of coherent sheaves on projective varieties after
being tensored with a sufficient large multiple of an ample line
bundle. In particular,  as a consequence of the result, we obtain
that the algebra $H^0(X,End (End E))$ coincides with $H^0(C,End
(End E|_C))$. This implies that the direct summands of $End E|_C$
are the restrictions to $C$ of the direct summands of $End E$.
Hence 3) follows if every saturated subsheaf $\Cal F$ of $End E$
with $\mu_H(\Cal F) = 0$ is a direct summand of $End E$. This last
statement is a consequence of 2). 2) implies that $\Cal F|_C$ is a
direct summand of $End E|_C$ and therefore, by the above,  $\Cal
F$ is a direct summand of $End E$. 4) follows from 3) and lemma
3.6. \hfill \ \qed
\enddemo

\proclaim{Corollary 3.10} Let $E$ be an absolute stable vector
bundle on a projective variety $X$.  Then $End
E=\bigoplus^l_{i=1}F_i$, where the $F_i$ are absolute stable
bundles with $\mu_H(F_i)=0$ for any $H$ in $P(X)$.
\endproclaim

\demo{Proof} The vector bundle $E$ is $H$-stable for all
polarizations $H$ of $X$. Fix a polarization $H$, lemma 3.9 4)
implies that $End E=\bigoplus^l_{i=1} F_i$ where all the $F_i$ are
$H$-stable with $\mu_H(F_i)=0$. Our claim is that the $F_i$ are
absolute stable vector bundles. \par Let $\Cal F$ be a coherent
subsheaf of one of the direct summands $F_i$ with $\text {rk} \Cal
F<\text {rk}F_i$. We need to show that $det \Cal F\in-K_{eff^+}$.
Lemma 3.8 almost gives the result, $det \Cal F\in-K_{eff}$.  If
$det \Cal F \notin -K_{eff^+}$ then $\mu_H(\Cal F)=0$. Lemma 3.9
3) implies that $\Cal F$ is a direct summand of $End E$ and hence
also of $F_i$. This is not possible since $F_i$ is $H$-stable.

\hfill \ \qed
\enddemo

\

\subhead 3.2 Holomorphic functions and flat bundles  \endsubhead

\

In section 2, we described a method to obtain holomorphic
functions on the universal cover $\rho:\tilde X \to X$ of a
complex manifold $X$. The method involved negative vector bundles
$E$ with a nontrivial cocycle $\alpha \in H^1(X,E)$ satisfying
$\rho^*\alpha=0$. In this subsection, we re-examine the method to
be able to apply it to vector bundles with very weak negativity
properties, see lemma 3.14 and proposition 3.16. One interesting
characteristic of these two results is that: if the vector bundle
$E$ satisfies the weak negativity  conditions described in the
results, the method fails to give nonconstant holomorphic
functions on $\tilde X$ only if $V$ is also a flat bundle. This is
interesting because Hodge theory and nonabelian Hodge theory
obtain holomorphic functions from flat bundles. Later in the
subsection 3.3, we will rely on this seemingly contradictory role
of flat bundles to describe the pullback map $\rho^*$. We will
 visit the production of holomorphic functions on the
universal covers of Kahler manifolds involving the existence of
flat bundles associated with infinite linear representations of
$\pi_1(X)$. We give  simple proofs for some  special cases. The
strongest result in this direction follows from [EKPR03] and is
described in observation 3.20.

 \

As previously announced, we now reexamine the method of producing
holomorphic functions developed in section 2. The goal is to be
able use the main idea of the method to get functions on $\tilde
X$ for the weakest possible "negativity" assumptions on the vector
bundle $V$.

\

\proclaim {Definition 3.11}  A sheaf $\Cal F$ on $X$ is
universally (generically) globally generated if the sheaf
$\rho^*\Cal F$ on the universal cover $\rho:\tilde X \to X$ is
(generically) globally generated.
\endproclaim

\

\proclaim {Lemma 3.12} Let $p: X' \to X$ be an infinite unramified
Galois covering of a complex manifold $X$ and $V$  a vector bundle
over $X$. If the kernel $p^*:H^1(X,V) \to H^1(X',p^*V)$ is
nontrivial then the vector bundle $p^*V$ on $X'$ has nonzero
sections.

\endproclaim

\input diagrams
\demo {Proof} Let $G$ be the Galois group of the covering and
$s\in H^1(X,V)$ be such that $p^*s=0$. The affine bundle
$A_{p^*s}$ is isomorphic to $p^*V$, but this isomorphism is not
$G$-equivariant with respect to the $G$-action on $p^*V$ whose
quotient is $V$. More precisely, there are two distinct actions of
$G$ on $p^*V$ which differ by affine transformations on $p^*V$.
One of the actions has as the quotient space the vector bundle $V$
and the other the affine bundle $A_s$. The action whose quotient
is $A_s$ can not preserve the zero section of $p^*V$ and hence
$p^*V$ has nontrivial sections. \hfill \ \qed
\enddemo

\

\proclaim {Lemma 3.13} Let $V$ be a vector bundle with a
nontrivial cocycle $\alpha \in H^1(X,V)$ such that
$\rho^*\alpha=0$. Then there is an universally globally generated
coherent subsheaf $\Cal F \subset V$ such that the cocycle
$\alpha$ comes from a cocycle $\beta \in H^1(X,\Cal F)$.
\endproclaim

\demo {Proof} It follows from lemma 3.12 that vector bundle
$\rho^*V$ has nontrivial sections. Let $\Cal F$ be the subsheaf of
$V$ whose stalk at $x\in X$ consists of the germs of the global
sections of $\rho^*V$ at one pre-image $\tilde x\in \rho^{-1}$.
Any choice of pre-image would give the same stalk since $\pi_1(X)$
acts on $\rho^*V$ and on $H^0(\tilde X,\rho^*V)$ as well.  The
sheaf $\Cal F$ is coherent because of the strong noetherian
property of coherent sheaves on complex manifolds. By construction
the sheaf $\Cal F$ is universally globally generated.\par

Let $i_*: H^1(X,\Cal F) \to H^1(X,V)$ and $q_*:H^1(X,V) \to
H^1(X,V/\Cal F)$ be the morphisms from the cohomology long exact
sequence associated with $0 \to \Cal F \to V \to V/\Cal F \to 0$.
The existence a cocycle $\beta \in H^1(X,\Cal F)$ with
$\alpha=i_*\beta$ follows if $q_*\alpha=0$. The extension $0 \to V
\to V_\alpha \to \Cal O \to 0$ associated with the cocycle
$\alpha$ induces the exact sequence:

 $$0 \to V/\Cal F \to V_\alpha/\Cal F \to \Cal O \to 0 \tag 3.1$$

 \noindent The
triviality of $q_*\alpha$ holds if (3.1) splits. The exact
sequence (3.1) is the quotient of the the exact sequence:

 $$0\to
\rho^*V/\rho^*\Cal F \to (\rho^*V)_{_{\rho^*\alpha}}/\rho^*\Cal F
\to \Cal O\to 0 \tag 3.2$$

\noindent via the action of $\Gamma=\pi_1(X)$ on
$(\rho^*V)_{_{\rho^*\alpha}}$ that gives $V_{_{\alpha}}$. The
extension of $\rho^*V$ associated with $\rho^*\alpha$ splits by
the hypothesis, but this splitting is not $\pi_1(X)$-invariant.
The splitting is given by a section  $s\in H^0(\tilde X,
\rho^*V_\alpha)$, that is not preserved by the $\pi_1(X)$-action.
On the other hand, this splitting  induces a $\Gamma$-invariant
splitting of (3.2) since $s-\gamma s\in H^0(\tilde X,\rho^*\Cal
F)$ and $\rho^*\Cal F/\Gamma=\Cal F$.
\enddemo

The next lemma is a flexible tool to produce holomorphic functions
on the universal coverings that will be a key ingredient of our
results.

\proclaim {Lemma 3.14} Let $\Cal F$ be a universally generically
globally generated coherent torsion free sheaf on a complex
manifold $X$ such that $\text { det} (\Cal F)^{-k}$ has a
nontrivial section. Then one of the following holds:\par

1) $\tilde X$ has a nonconstant holomorphic function.\par

2) $\Cal F$ is the sheaf of sections of a flat bundle, $\Cal F
\cong \Cal O(\tau)$.
\endproclaim

\demo {Proof} Let $s_1$,...,$s_r$ be a collection of sections of
$\rho^*\Cal F$ generating $\rho^*\Cal F$ generically, where $r$ is
the rank of $\Cal F$. From the sections $s_1$,...,$s_r$ one gets a
nontrivial section of $\text {det}\rho^*\Cal F$
$s=s_1\wedge...\wedge s_r\in H^0(\tilde X,\text {det}(\rho^*\Cal
F))$, the pairing of $s^{\otimes n}$ with a nontrivial section $t
\in H^0(\tilde X,\text { det}(\rho^*\Cal F)^{-n})$ gives a
holomorphic function $f$ on $\tilde X$. By hypothesis the function
$f$ is nonzero on a open set of $\tilde X$. The function $f(p)=0$
at $p\in \tilde X$ in case: i) $s_1(p)$,...,$s_r(p)$ do not
generate $\rho^*{\Cal F_p}/m_p\rho^*{\Cal F_p}$ or ii) $t(p)$ is
zero. \par Suppose statement 1) does not hold. Then $f$ must be a
nonzero constant function, which implies that
$s_1(p)$,...,$s_r(p)$ are linear independent at all $p\in \tilde
X$. Hence the morphism $(s_1,...,s_r):\Cal O^r \to \Cal \rho^*\Cal
F$ induced by the sections is an isomorphism. The nonexistence of
holomorphic functions on $\tilde X$ implies that all sections of
$\rho^*\Cal F=\Cal O^r$ are constant. The linear action of
$\pi_1(X)$ on $H^0(\tilde X,\Cal O^r)=\Bbb C^r$ gives a
representation $\tau: \pi_1(X) \to GL(r,\Bbb C)$ and $\Cal F $ is
the sheaf of sections of the flat vector bundle $\tilde X
\times_\tau \Bbb C^r$.

 \hfill \ \qed
\enddemo

\proclaim{Corollary 3.15} Let $X$ be a projective variety such
that $H^0(\tilde X,\Cal O)=\Bbb C$. If $E$ is a vector bundle such
that $det E^* \in K_{eff}$ then $E$ is not universally generically
global generated unless $det E$ has finite order in $Pic X$ and
$E$ is flat.
\endproclaim

The next result is an application of lemma 3.14 for vector
bundles.

\proclaim {Proposition 3.16} Let $E$ be an absolutely stable
vector bundle over a projective manifold $X$ with $det E=\Cal O$.
If there is a nontrivial cocycle $\alpha \in H^1(X,E)$ such that
$\rho^*\alpha=0$ then one of the following possibilities
holds:\par

1) $\tilde X$ has nonconstant holomorphic functions. \par

2) $E\cong \Cal O(\tau)$ and $\alpha$ is contained in the image of
$H^1(X,\Bbb C(\tau))$ in $H^1(X,E)$.
\endproclaim

\demo {Proof} Lemma 3.13 states that there is a nontrivial
universally globally generated coherent subsheaf $\Cal F$ of $E$
such that $\alpha$ is contained in the image of $H^1(X,\Cal F)$ in
$H^1(X,E)$. If $\text {rk}\Cal F< \text {rk} E$ then the absolute
stability of $E$ would imply that  the line bundle $(det \Cal F)^*
\in K_{eff^+}$. Hence $\Cal F$ is not a flat vector bundle and by
lemma 3.14 $\tilde X$ must have nonconstant holomorphic functions.

If $\text {rk}\Cal F= \text {rk} E$ then $E$ is universally
generically globally generated. Since $\tilde X$ has no
nonconstant holomorphic functions, lemma 3.14 gives that $E$ is a
flat vector bundle $\tilde X \times_\tau \Bbb C^r$, for some
representation $\tau: \pi_1(X) \to GL(r,\Bbb C)$. Notice that for
$\alpha\in H^1(X,E)$ with $\rho^* \alpha= 0$ the section $s\in
H^0(\tilde X, \Cal O^r)$ ($\rho^* E \cong \Cal O^r$) with $d s =
\rho^*\alpha$ is constant. Hence we obtain that $s$ belongs to the
image of $s' \in H^1(X,\Bbb C(\tau))$ under a natural map
$H^1(X,\Bbb C(\tau))\to H^1(X,E)$. \hfill \ \qed
\enddemo

\

Both lemma 3.14 and proposition 3.16 give a method to obtain
holomorphic functions on $\tilde X$ from vector bundles that are
not flat. The results that follow give a description of how flat
bundles can produce holomorphic functions on the universal cover.
Recall the notation described in the introduction to this section,
a representation $\tau : \pi_1(X) \to GL(m,\Bbb C)$ defines the
flat vector bundle $\Cal O(\tau)$ on $X$. We denote the sheaf of
sections of $\Cal O(\tau)$ also by $ \Cal O(\tau)$ and the
subsheaf of locally constant sections by $\Bbb C(\tau)$. The
imbedding $\Bbb C(\tau) \hookrightarrow \Cal O(\tau)$ induces a
map of the cohomology groups.

\proclaim{Proposition 3.17} Let $X$ be a complex manifold and
$\tau : \pi_1(X) \to GL(m,\Bbb C)$ be a representation of
$\pi_1(X)$. If $\tilde X$ has no nonconstant holomorphic functions
then the map $H^1(X,\Bbb C(\tau)) \to H^1(X,\Cal O(\tau))$ is an
imbedding.
\endproclaim

\demo{Proof} Consider the exact sequence of sheaves on $X$
associated with the differential $d$:  \diagram [size=1.7em]
0&\rTo&\Bbb C(\tau)&\rTo&\Cal O(\tau)&\rTo^{_d}&d\Cal
O(\tau)&\rTo&0
\enddiagram

\noindent  where $d\Cal O(\tau)$ is the image  subsheaf in
$\Omega^1(X)\otimes \Cal O(\tau)$. From the long cohomology exact
sequence, we have:

$$H^0(X,d\Cal O(\tau)) \to H^1(X,\Bbb C(\tau))\to H^1(X,\Cal O(\tau)) $$

If the second map has a nontrivial kernel, then $H^0(X,d\Cal
O(\tau))$ is nonzero. Any section of\par \noindent $H^0(X,d\Cal
O(\tau))$ induces a closed holomorphic (1,0)-form on $\tilde X$
with  values in $\Cal O^n$ and by integration a set of nonconstant
holomorphic functions. \hfill \ \qed
\enddemo

\

\proclaim{Corollary 3.18} Let $X$ be a complex manifold and $\tau
: \pi_1(X) \to GL(m,\Bbb C)$ be a representation of $\pi_1(X)$. If
$h^1(X,\Bbb C(\tau))> h^1(X,\Cal O(\tau))$ then $\tilde X$ has
nonconstant holomorphic functions.
\endproclaim

\demo{Proof} Suppose the corollary did not hold then  $H^0(\tilde
X,\Cal O)=\Bbb C$. Hence $H^0(X,\Bbb C(\tau))=H^0(X,\Cal O(\tau))$
 since any holomorphic section of $\rho^*\Cal O(\tau)$ on $\tilde X$
 is constant. Apply proposition 3.17 to get the contradiction. \hfill
\ \qed
\enddemo

\

\proclaim{Proposition  3.19} Let $X$ be a complex manifold such
that $\tilde X$ has no nonconstant holomorphic functions. Then the
following properties hold:\par
  1) For any bundle $E$ there is at most one representation $\tau$, up to
conjugation, such that $E = \Cal O(\tau)$.\par
  2) If $X$ is Kahler then $H^1(X,\Bbb C(\tau)) = H^1(X,\Cal O(\tau))=0 $ for any unitary representation
 $\tau $ of $\pi_1(X)$. In particular,  $H^1(X,\Cal O) = 0$, $Pic_{
0} X\otimes \Bbb Q = 0$ and  $Pic X\otimes \Bbb Q = NS(X)\otimes
\Bbb Q$.
\endproclaim

\demo{Proof} 1) The structure of a bundle  $\Cal O(\tau)$ on $E$
is the same as a flat connection on $E$. Suppose there were two
different structures $\Cal O(\tau)$ and $\Cal O(\tau')$ on $E$,
they would induce two flat (1,0)-connections whose difference is
an non-zero element of $H^0(X, dEnd \Cal O(\tau))$. The desired
conclusion follows from the argument in proof of proposition 3.17.
\par

2) Notice that for a unitary representation $\tau$ the cohomology
of $\Bbb C(\tau)$ and $\Cal O(\tau)$ satisfy the Hodge
decomposition. In particular, there is an isomorphism of vector
spaces $H^1(X,\Bbb C(\tau))= H^1(X,\Cal O(\tau)) \oplus
H^0(X,\Omega^1\otimes \Cal O(\bar \tau))$. Corollary 3.18 implies
that $H^0( X,\Omega^1\otimes \Cal O(\bar \tau)) = 0$ and by the
Hodge conjugation isomorphism, it follows that $H^1(X,\Cal
O(\tau))=0$. As a special case, we obtain $H^1(X,\Cal O) = 0$.

 \hfill
\ \qed
\enddemo

\proclaim {Observation 3.20}  The unpublished paper [EKPR03] has
implicit the following consequence: if a smooth Kahler variety $X$
has an infinite linear representation of the fundamental group
then its universal cover has  nonconstant holomorphic functions.
For projective surfaces this result appears in [Ka97]. For Kahler
manifolds, our three last results follow this consequence. On the
other hand, the proof of our results is more direct and
significantly simpler.
\endproclaim

\

\noindent Remark: The method of [EKPR03] using non-abelian Hodge
theory to construct holomorphic functions on the universal covers
$\tilde X$ requires that the base manifold $X$ is Kahler. If $X$
is not Kahler the properties of the fundamental group of $X$ can
not guarantee the existence of non-constant holomorphic functions
on $\tilde X$. This follows from the results of Taubes on
anti-self-dual structures on real 4-manifolds [Ta92]. Taubes
showed that every finitely presented group is the fundamental
group of a compact complex 3-fold $X$ that has a foliation by
$\Bbb P^1$ with normal bundle $\Cal O(1)\oplus \Cal O(1)$. This in
turn implies that the universal cover $\tilde X$ has no
non-constant holomorphic functions. The universal cover $\tilde X$
also has a foliation by $\Bbb P^1$ with normal bundle $\Cal
O(1)\oplus \Cal O(1)$. Any one of these $\Bbb P^1$ has a 2-concave
and hence pseudoconcave neighborhood since their normal bundle is
Griffiths-positive [Sc73]. The conclusion follows since a complex
manifold with a pseudoconcave open subset has only constant
holomorphic functions. Moreover the variety $\tilde X$ with $X$
being a twistor space for a sufficiently generic anti-self-dual
metric on the underlying $4$-dimensional variety has no
meromorpihic functions. Indeed the field of meromorphic functions
on $\tilde X$ is always a subfield of the field of meromorphic
functions in the normal neighborhood of $\Bbb P^1$ and the latter
is always a subfield of $C(x,y)$ and consists of constant
functions only for a generic neighborhood of $\Bbb P^1$ with
normal bundle $\Cal O(1)\oplus \Cal O(1)$.

\

\subhead 3.3  Pullback map for line bundles \endsubhead

\

We describe the  implications of the absence of nonconstant
holomorphic functions on the universal cover $\tilde X$ on the
pullback map for line bundles $\rho^*: \text{Pic} (X) \to
\text{Pic} (\tilde X)$.

\proclaim {Definition 3.21}  The cone of divisors on $X$ generated
by the divisors which become effective on $\tilde X$ is denoted by
$\tilde K_{eff}$. Following 3.2, $\tilde K_{eff^+}$ is also
defined.
\endproclaim

If $\tilde X$ has no non-constant holomorphic functions then the
cone $\tilde K_{eff}$ contains $K_{eff}$ but does not contain any
elements from $- K_{eff^+}$.  Suppose  $- K_{eff^+} \cap \tilde
K_{eff}\neq \emptyset$ then there is an divisor effective $D$ of
$X$ such that both line bundles $\rho^* \Cal O(D)$ and $\rho^*\Cal
O(-D)$ have nontrivial sections. The pairing of these sections
gives a non-constant holomorphic function. In particular, the
image of $\rho^*: \text{Pic} (X) \to \text{Pic} (\tilde X)$ is
nontrivial and there are the following possibilities:

 \

 P1. The cone $\tilde K_{eff^+}$ is separated by a hyperplane
 $L$ from $ - K_{eff^+}$.

 P1'. The cone $\tilde K_{eff}$ coincides with $K_{eff}$.
  This is a special case of P1. In particular, it holds if
  $H^0(\tilde X,\Cal O)=\Bbb C$ and $\text {Pic} (X) = \Bbb Z$.

 P2. The closure of the cone $\tilde K_{eff}$ in $\text {Pic} (X)\otimes \Bbb R$
  intersects with the closure of $- K_{eff}$ outside of $0$.

  \

The following result  describes the kernel of $\rho^*: \text{Pic}
(X) \to \text{Pic} (\tilde X)$ for Kahler manifolds whose
universal cover has no nonconstant holomorphic functions.

\proclaim{Proposition 3.22} Let $X$ be a Kahler manifold such that
$H^0(\tilde X,\Cal O)=\Bbb C$. Then the kernel of the pullback map
$\rho^*: Pic X \to Pic \tilde X$ is finite and its elements
correspond to flat bundles associated with finite characters.
\endproclaim

\demo{Proof} Let $L$ be a line bundle in the kernel of $\rho^*$
and let $i :\rho^* L \to \Cal O_{\tilde X}$ be an isomorphism with
the trivial line bundle. The isomorphism $i$ is not equivariant
with respect to the natural $\pi_1(X)$-actions on $\Cal O_{\tilde
X}$ and on $\rho^* L$ giving respectively $\Cal O_X$ and $L$ on
$X$. Hence there is a $g\in \pi_1(X)$ such that the map
$(gi)i^{-1}\neq \text {Id}:\Cal O_{\tilde X} \to \Cal O_{\tilde
X}$. If $(gi)i^{-1}\neq \text {cId}$ for some constant $c$ then we
obtained a nonconstant holomorphic function on $\tilde X$, which
can not happen.\par Therefore, we have an association of elements
of $\pi_1(X)$ with nonzero constants. This association defines a
representation $\pi_1(X)\to \Bbb C^*$ and this representation has
to be finite, since $H^0(\tilde X, \Cal O)=\Bbb C$ implies that
$H^1(X,\Bbb C)$ vanishes. The line bundles on the kernel of
$\rho^*$ are uniquely determined by the representations described
above. Thus, $\text {Ker}(\rho^*)$ is dual to $\pi_1(X)^{ab}$
which is a finite group.\hfill \ \qed
\enddemo

\

\subhead 3.4  Pullback map for vector bundles \endsubhead

\

Assuming that $\tilde X$ has no nonconstant holomorphic functions,
we use the previous results to describe the pullback map on the
moduli spaces of absolutely stable vector bundles on $X$. Our
results are mostly for the spaces of absolutely stable bundles but
they have a generalization for the spaces of $H$-stable bundles,
if extra conditions on $\tilde K_{eff}$ are added.

In order to describe the local behavior of the pullback map, it
necessary to recall some facts from the theory of deformations of
a given vector bundle $E$ on $X$. The deformation of a vector
bundle $E$ over an arbitrary variety splits into the deformation
of the projective bundle $\Bbb P(E)$ plus a deformation of the
line bundle $\Cal O_{\Bbb P(E)}(1)$ over $\Bbb P(E)$. The
deformations of $\Bbb P(E)$ in the case of a smooth $X$ are
parameterized by an analytic subset $B_E\subset H^1(X,End E), 0\in
B_E$ with the action of the group of  relative analytic
automorphisms $Aut (E)$ of the bundle $\Bbb P(E)$ on $B_E$. The
latter is induced from the natural linear action $Aut (E)$ on
$H^1(X,End E)$(with adjoint fiberwise action of $PGL(n)$ on the
fiber of the bundle $End_E$). Thus, non-isomorphic bundles (with
respect to identical automorphism on $X$) in the local
neighborhood of $E$ are parameterized by the orbits of the group
$Aut E$ with Lie algebra $H^0(X, End_0(E) $ in $H^1(X,End E)$.

The space $H^1(X,End E)$ plays a role of the formal tangent space
$T_0(B_E)$ at the point $0\in B_E$. Natural splitting $End E =
End_0E \oplus \Cal O$ induces a splitting $H^1 (X,End E) =
H^1(X,End_0 E) \oplus H^1(X,\Cal O)$. The local deformation scheme
of $E$ maps onto a local deformation scheme of $\Bbb P(E)$ with a
fiber which is locally isomorphic to  $H^1(X,\Cal O)$. $H^1(X,\Cal
O)$ parameterizes the (non-obstructed) deformation scheme of line
bundles $\Cal O_{\Bbb P(E)}(1)$ in $\text {Pic} (\Bbb P(E))$ over
the deformation scheme of $\Bbb P(E)$ which is generically
obstructed.

\

Let $p:\Cal E \to \Delta$ be an analytic family, over the disc
$\Delta$, of vector bundles on $X$ with $E_t=p^{-1}(t)$ as its
members. The family $E_t$ gives a deformation of $E=E_0$ and has
associated with it a 1st-order deformation cocycle $s\in H^1(X,End
E)$.

\proclaim{Lemma 3.23} Let $p:\Cal E \to \Delta$ be family  of
vector bundles on $X$ that is nontrivial at $t=0$. If the pullback
family $\tilde p:\rho^*\Cal E \to \Delta$ is locally trivial then
the kernel of $\rho^*: H^1(X,End E) \to H^1(\tilde X,End \rho^*E)$
is nontrivial.
\endproclaim

\demo{Proof } The 1st-order deformation  cocycle $s\in H^1(X,End
E)$ associated with the family $E_t$ is nontrivial since the
family $E_t$ is nontrivial at $t=0$. The nontrivial cocycle $s$ is
in the kernel of $\rho^*$ since $\rho^*s$, the 1st-order
deformation cocycle associated with  the locally trivial family
$\tilde p:\rho^*\Cal E \to \Delta$, is trivial. \hfill \ \qed
\enddemo

\

\proclaim{Lemma 3.24}  Let $p:X' \to X$ be an unramified Galois
covering of a smooth projective manifold $X$ and $E$ a vector
bundle on $X$. Then $H^0(\tilde X,End_0 p^*E) \neq 0$ if one of
the following holds:\par 1) The kernel of $p^* :H^1(X,End_0E)\to
H^1(X',End_0 p^*E)$ is nontrivial.\par 2) $H^0(X',\Cal O)=\Bbb C$
and there is a pair of vector bundles $E$ and $F$ such that $p^*
F=p^* E$ but $F\neq E\otimes \Cal O(\chi)$ for any character $\chi
: \pi_1(X)\to \Bbb C^*$.
\endproclaim

 \demo{Proof} Assume that 1) holds then  $H^0(\tilde X,End_0
 p^*E)\neq 0$ follows from lemma 3.12 ($p$ must be an infinite
 unramified covering of $X$ by Lemma 2.2).\par

 If 2) holds then there is an isomorphism $i:p^* E \to p^* F$ and
 $F\neq E\otimes \Cal O(\chi)$ for any character
$\chi : G\to \Bbb C^*$. Let $G$ be the Galois group of the
covering. The isomorphism $i$ is not $G$-equivariant since
otherwise it would descent to an isomorphism $i':E \to F$ on $X$.
Consider the two possible cases: 1) there is a $g\in G$ such that
$g^{-1}i^{-1}gi : p^*E \to p^*E$ is a non-scalar endomorphism.
Then $g^{-1}i^{-1}gi$ is a nontrivial element in
$H^0(X',End_0p^*E)$. 2) For all $g\in G$ the endomorphism
$g^{-1}i^{-1}gi$ of $p^*E$ is scalar. Since $X'$ has no
nonconstant holomorphic functions, the following holds:
 $g^{-1}i^{-1}gi=\chi(g)\text {Id}$, $\chi(g) \in \Bbb C^*$. Therefore, the
 map $\chi: G \to \Bbb C^*$ defines a character of $G$ and
  $F = E\otimes \Cal O(\chi)$ which can not happen, since it contradicts the assumption.\hfill \ \qed
\enddemo

\

Let $\rho^*_0:Mod_0(X) \to Vect(\tilde X)$ be the pullback map,
where $Mod_0(X)$ is the  moduli space of absolutely stable vector
bundles on $X$.  We denote points in $Mod_0(X)$ by the same
letters as the corresponding vector bundles.

\proclaim {Proposition 3.25} Let $X$ be a projective manifold such
that its universal cover $\tilde X$ has no nonconstant holomorphic
functions. If $E$ is an absolutely stable vector bundle on $X$
satisfying $H^0(\tilde X,End_0\rho^* E) = 0$, then:\par a) The
pullback map $\rho_0^*:Mod_0(X) \to Vect(\tilde X)$ is a local
embedding at $E$.\par b) For any  absolutely stable bundle $E$
there are only finite number of bundles $F$ with $\rho^*E = \rho^*
F$ and  $E = F\otimes \Cal O(\chi)$ with $\chi$ a character of
$\pi_1(X)$.\par
\endproclaim

\demo{Proof} To prove part a) it is enough to show that the
tangent map to $\rho_0^*$ at $E$,  $\rho^*:H^1(X, End E) \to
H^1(\tilde X, End \rho^*E)$ is injective. The injectivity of
$\rho^*:H^1(X, End E) \to H^1(\tilde X, End \rho^*E)$ follows from
$H^0(\tilde X,End_0\rho^* E) = 0$ and  lemma 3.24 1).\par

Lemma 3.24 2) and the finiteness of the character group of
$\pi_1(X)$ imply  part b). The finiteness of the character group
follows from $H^0(\tilde X,\Cal O)= \Bbb C$.

\hfill \ \qed
\enddemo

To conclude, we consider the case when  $E$ is an absolute stable
but $H^0(\tilde X,End_0 E) \neq 0$.

\proclaim{Theorem 3.26}  Let $X$ be a projective manifold such
that its universal cover $\tilde X$ has no nonconstant holomorphic
functions. If $E$ is an absolute stable vector bundle on $X$
satisfying $H^0(\tilde X,End_0 \rho^*E) \neq 0$ then associated to
$E$ is a normal subgroup $\pi_1(X')\subset \pi_1(X)$ corresponding
to finite unramified covering $p :X' \to X$ with universal cover
$\rho': \tilde X \to X'$ satisfying:

\par
i) $p^*E\simeq E'_1\otimes \Cal O(\tau_1)\oplus ...\oplus
E'_m\otimes \Cal O(\tau_m)$, $\tau_i:\pi_1(X') \to GL(k,\Bbb C)$
with $k\ge 1$.

\par ii) $H^0(\tilde X, End_0{\rho'}^*E'_i)=0$ for
all $i=1,...,m$.

\par iii) The natural action of the finite group $G = \pi_1(X)/\pi_1(X')$
 on $X'$ extends to the action on $p^*E$ which permutes subbundles
$ E'_i\otimes \Cal O(\tau_i)$ and this action gives the imbedding
 of $G$ into $S_m$.

\par iv) Let $G_1\subset G$ be subgroup which acts
identically on $E'_1\otimes \Cal O(\tau_1)\subset p^*E$ then
$E'_1\otimes \Cal O(\tau_1)$ descends to the bundle $E'_1\otimes
\Cal O(\tau_1')$ on $X_1 = X'/G_1$
 with $p_1 : X_1\to X$ being a nonramified covering
 of degree $rk E / (rk E'_1\otimes \Cal O(\tau_1))$
 and $E = p_* E'_1\otimes \Cal O(\tau_1')$

\endproclaim

\demo{Proof} Consider the subsheaf $\Cal A$ of $End \rho^*E$
generated by the global sections of $End \rho^*E$. The sheaf $End
\rho^*E$ is a sheaf of matrix algebras and $\Cal A$ is a sheaf of
subalgebras since we can add and multiply sections.  We claim
that: $A=H^0(\tilde X,\Cal A)=H^0(\tilde X, End \rho^*E)$ is
finite dimensional and isomorphic to a sum of $m$ copies of the
algebra of $k\times k$ matrices $M(k)$ for a $k<r=\text {rk } E$;
the action of $\pi_1(X)$  on the  algebra $A$ has no nontrivial
$\pi_1(X)$ invariant ideals ($\pi_1(X)$ acts transitively on the
$m$ direct summands $M(k)\subset A=\bigoplus_{i=1}^mM(k)$). We
also claim that the  action of the algebra $A$ on $\rho^*E$ is
such that each direct summand of $A$ acts on each fiber
$(\rho^*E)_x\simeq \Bbb C^r$, $x \in \tilde X$, as the same
multiple $lM(k)$ of the standard representation of $M(k)$,
$r=lmk$.
\par

As mentioned above, the $\pi_1(X)$-action on $A$ permutes the m
simple direct summands $M(k)$ and therefore gives a homomorphism
$\sigma:\pi_1(X) \to S_m$. Associated with the normal subgroup
$Ker(\sigma)$ is a finite unramified Galois covering of $X$, $p:X'
\to X$, with $\pi_1(X')=Ker(\sigma)$. By construction the direct
summands of $A$ are $\pi_1(X')$-invariant. Thus, $\rho^*E=\tilde
E_1\otimes \Cal O^k \oplus ... \oplus {\tilde E}_m\otimes \Cal
O^k$, where $\tilde E_i={\rho'}^*E'_i$, $E'_i$ is a vector bundle
of rank $l$ on $X'$ and $\rho':\tilde X \to X'$ is the universal
cover of $X'$. On $X'$ the bundle $p^*E$ decomposes into $p^*E
=E'_1\otimes \Cal O(\tau_1)\oplus ... \oplus E'_m\otimes \Cal
O(\tau_m)$ ($\tau_i:\pi_1(X') \to GL(k,\Bbb C)$) giving i). Note
that the group $G=\pi_1(X)/\pi_1(X')$ acts on $p^*E$ permuting
transitively the direct summands thus proving iii). Part ii)
follows from $H^0(\tilde X, End({\rho'}^*E'_i\otimes \Cal
O^k))=M(k)$ since if $H^0(\tilde X, End_0{\rho'}^*E'_i)\neq 0$ the
group of global sections $H^0(\tilde X, End({\rho'}^*E'_i\otimes
\Cal O^k))$ would be larger. Consider the group $G_1\subset G$
which stabilizes $E'_1\otimes \Cal O(\tau_1)$. Then $E'_1\otimes
\Cal O(\tau_1)$ descends to $X_1/G_1$. The bundle $p^*E$ also
descends to $X_1$ and it decomposes into a direct sum $E'_1\otimes
\Cal O(\tau_1) + E_N$. Let $p_1 : X_1 \to X$ be a corresponding
covering of $X$. Consider the direct image $p_{1*}E'_1\otimes \Cal
O(\tau_1) $ on $X$. We want to show that $p_{1*}E'_1\otimes \Cal
O(\tau_1) = E$. Consider also $p_{1,*}p^* E $ which has a natural
decomposition as $E + E_c$. Natural projection $i_1^*: p^* E \to
E'_1\otimes \Cal O(\tau_1)$ which is identity on $E'_1\otimes \Cal
O(\tau_1)$ induces a map $p_{1,*} p^* E \to p_{1,*}E'_1\otimes
\Cal O(\tau_1)$. Denote by $R$ the restriction of $R$ on the
direct summand $E\in p_{1,*}p^* E $. We want to show that $R$ is
an isomorphism. It follows from the fiberwise description of $R$.
Let $C^{mlk}_x = \Sigma C^{mk}_i$ be a direct decomposition of the
fiber of $p^*E$ at $x\in X'$ into the sum of the fibers of the
direct summands $E'_i\otimes \Cal O(\tau_i)$ and $g_i\in G$ be the
representatives of cosets $G/G_1$. Then for a $x'\in X$ its
pre-image $p_1^{-1}x' \subset X_1$ is equal to $\bigcup g_ix_1$
and the fiber of $ E'_1\otimes \Cal O(\tau_1) $ over $g_ix_1$ is
naturally isomorphic to $C^{mk}_i$. Now the map $R$ becomes the
trace map for the action of $G$ on $p^*E$ which implies that $R$
is fiberwise isomorphism. This proves iv).

\

\noindent Claim:  $A=H^0(\tilde X,\Cal A)$ is a subalgebra of the
matrix algebra $M(r)$, $r= \text {rk} E$.\smallskip

First, we prove the finite dimensionality of $A$. The sheaf $\Cal
A$ is invariant under the action of $\pi_1(X)$ and defines a
coherent subsheaf $\Cal A'\subset End E$ on $X$ with $\Cal
A=\rho^*\Cal A'$. The absolute stability of $E$ implies that $det
\Cal A' \in -K_{eff}$ by corollary 3.8. The lemma 3.14 implies
that $\Cal A'$ is isomorphic to the sheaf of sections of a flat
vector bundle since $(det \Cal A')^{-k}$ has a nontrivial section
for some $k$ but $H^0(\tilde X,\Cal O)=\Bbb C$. Hence $\Cal
A=\rho^*\Cal A'=\Cal O^q$.  It follows from $H^0(\tilde X,\Cal
O)=\Bbb C$, that the algebra $A$ is finite dimensional.
\par

We want to show that  the algebra $A\cong \Cal A\otimes
k(x)\subset End \rho^*E\otimes k(x) \cong M(r)$, where $x$ is any
point in $\tilde X$ and $k(x)$ is the residue field at $x$.
Consider the exact sequence $0\to \Cal A\otimes \Cal I(x) \to \Cal
A \to \Cal A\otimes k(x) \to 0$, $\Cal I(x)$ the ideal sheaf of
the point $x$. Since $\Cal A$ is globally generated by its global
sections it follows that the morphism $A=H^0(\tilde X,\Cal A) \to
H^0(X,\Cal A\otimes k(x))=\Cal A\otimes k(x)$ is a surjection. If
the morphism is also an injection we get the desired isomorphism
$A\cong \Cal A\otimes k(x)$. The injectivity follows from
$H^0(\tilde X,\Cal A\otimes \Cal I(x))=0$, which holds since the
argument of the previous paragraph implies that any nontrivial
section $s$ of $End\rho^*E$ is nowhere vanishing.

\

\noindent Claim: The algebra $A$ is semisimple.\smallskip The
semisimplicity of $A$ is equivalent to the maximal nilpotent ideal
$I_m$ of $A$ being the zero ideal. The algebra $A$ comes with a
natural $\pi_1(X)$-action. The maximal nilpotent ideal is a
$\pi_1(X)$-invariant ideal of $A$. Every nontrivial
$\pi_1(X)$-invariant ideal $I$ of $A$ defines naturally a
nontrivial subsheaf $\Cal I' \subset \Cal A' \subset End E$.
Suppose that the ideal $I_m$ is nontrivial and consider the
subsheaf $\Cal I'_mE$ of $E$. The nilpotent condition, $I^k_m=0$
for some $k$, implies that $\text {rk}(\Cal I'_mE)<\text {rk}
E$.\par

 Consider the exact sequence $0 \to \Cal K \to \Cal
I_m'\otimes E \to \Cal I_m'E \to 0$. We will show that the
subsheaf $\Cal K$ is an $H$-destabilizing subsheaf of $\Cal
I_m'\otimes E$ for any polarization of $X$. This can not be, since
$E$ and $\Cal I_m'$ are $H$-semistable ($\Cal I_m'$ is a flat
bundle) and hence $\Cal I_m'\otimes E$ is also $H$-semistable (the
tensor product of two semistable sheaves is semistable). We need
to get the destabilizing inequality $[\text {rk}(\Cal I_m'\otimes
E) det \Cal K-\text {rk}\Cal K det (\Cal I_m'\otimes
E)].H^{n-1}>0$. Using $det \Cal K=\text {rk} \Cal I_m'det E-det
\Cal I'E$ and $\text {rk}Edet \Cal I_m'E-\text {rk}\Cal
I_m'EdetE\in K_{eff^+}$ ($E$ is an absolute stable bundle), it
follows that $\text {rk}(\Cal I_m'\otimes E) det \Cal K-\text
{rk}\Cal K det (\Cal I_m'\otimes E)\in -K_{eff^+}$. Hence we
obtain the desired contradiction, which implies that the
$\pi_1(X)$-invariant ideal $I_m$ must be the zero ideal and $A$ is
semisimple.

\

\

\noindent Claim: $A$ has no proper $\pi_1(X)$-invariant
ideals.\smallskip

We proved that $A$ is semisimple and hence $A =\Sigma^m_{i=1}
M(r_i)$ with $r_1+...+r_m\le r$. First, we note that if $I$ is an
ideal of $A=\Sigma^l_{i=1} M(n_i)$ ($A$ acts on $\Bbb C^r$) such
that $I\Bbb C^r=\Bbb C^r$ then $I=A$. We prove the claim by
showing that for any nonzero $\pi_1(X)$-invariant ideal $I$ of $A$
the equality $I=A$ or equivalently $\Cal I'E=E$ must hold. If
$\text {rk}\Cal I'E=\text {rk}E$ then at some $x\in \tilde X$,
$I_x=I\otimes k(x)\subset End(\rho^*E)_x$ is such that
$I_x(\rho^*E)_x=(\rho^*E)_x$ and hence $I=A$. If $I$ is such that
$\text {rk}\Cal I'E<\text {rk}E$ then the argument of the
paragraph above proves that $I=0$ and hence also the claim.

\

\noindent Claim: The algebra $A$ is equal to $mM(k)$ and the
representation of each $M(k)$ is a multiple of a standard rank $k$
representation of $M(k)$.\smallskip

The algebra $A$ is as noted before $A=\Sigma^m_{i=1} M(r_i)$ with
$r_1+...+r_m\le r$. Since each $M(r_i)$ is simple it follows that
the action of $\pi_1(X)$ preserves the ideals of $A$ corresponding
to the sums of all the $M(r_i)$ with $r_i$ equal to a fixed $k$.
Therefore all the $r_i$ are equal to the same $k$ since any
$\pi_1(X)$-invariant ideal $I$ of $A$ is either trivial or the
full $A$.\par

Finally, we show that  the representation of each $A=M(k)$ in
$M(n)$ is a multiple of the standard representation. Any
irreducible representation of $M(k)$ is the standard
representation or the zero representation. The presence of a zero
representation as an irreducible component of the representation
of $M(k)$ in $M(n)$ would imply that $\Cal A'E\neq E$ which is not
possible from the discussion above.\par

 \hfill \ \qed
\enddemo

\

The following lemma that follows from our results and observation
3.20.

\proclaim {Lemma 3.27} If $\tilde X$ has no nonconstant
holomorphic functions then $H^1(X,\Bbb C(\tau))=H^1(X,\Cal
O(\tau))=0$ for any linear representation $\tau$ of $\pi_1(X)$.
\endproclaim
\demo {Proof}  Let $f:X' \to X$ of $X$ be the covering
corresponding to the kernel $G \subset \pi_1(X)$  of the
representation   $\tau$. The hypothesis $H^0(\tilde X,\Cal O)=\Bbb
C$ and observation 3.20 imply that $\tau$ is finite. It also
follows from $H^0(\tilde X,\Cal O)=\Bbb C$ that $H^1(X' ,\Bbb
C^k)=H^1(X' ,\Cal O^k)=0$. Since $\tau$ is finite  the covering
$f$ is also finite  and hence the conclusion follows from the
imbeddings $f^*:H^1(X, \Bbb C(\tau)) \to H^1(X' ,\Bbb C^k)$ and
$f^*:H^1(X, \Cal O(\tau)) \to H^1(X' ,\Cal O^k)$ . \hfill \ \qed
\enddemo

\

\proclaim {Theorem B} Let $X$ be a projective manifold whose
universal cover has only constant holomorphic functions. Then:
\par

a) The pullback map $\rho_0^*:Mod_0(X) \to Vect(\tilde X)$  is a
local embedding.
\par b) For any absolutely stable bundle $E$ there are
only finite number of bundles $F$ with $\rho^*E = \rho^* F$.\par
c) Moreover, there is a finite unramified cover $p: X'\to X$
associated with $E$ of degree $d\leq rk E !$ with universal
covering ${\rho'}:\tilde X \to X'$. On $X'$ there is  a collection
of vector bundles $\{E'_i\}_{i=1,...,m}$ on $X'$ with $H^0(\tilde
X,End_0{\rho'}^*E'_i)=0$ such that $\rho^*F\simeq \rho^*E$ if and
only if:
 $$p^*F =E'_1\otimes \Cal O(\tau_1)\oplus ... \oplus E'_1\otimes \Cal
 O(\tau_m)$$
 The bundles $\Cal O(\tau_i)$  are flat bundles associated with finite linear
representations of $\pi_1(X')$ of a fixed rank $k$ with $rk E|k$.
\endproclaim

\demo {Proof} The case for vector bundles $E$ such that
$H^0(\tilde X,End_0\rho^* E) = 0$ was done in proposition 3.25. We
proceed to consider the case $H^0(\tilde X,End_0\rho^* E)\neq
0$.\par a) Let $Mod_0(X,V)$ be the moduli space of absolutely
stable bundles with the same Chern classes as $V$. The formal
tangent space of $Mod_0(X,V)$ at $V$ is given by $H^1(X,EndV)$.
The vector bundle $EndV$ is semistable with $det EndV=\Cal O$ and
is the direct sum $End E=\bigoplus_{i=1}^lF_i$ of absolutely
stable bundles with $\mu_H(F_i)=0$ by corollary 3.10. The kernel
of the tangent map $\rho^*:H^1(X,EndV)\to H^1(\tilde
X,End\rho^*V)$ is the direct sum of the kernels of $\rho^*_i\equiv
\rho^*:H^1(X,F_i)\to H^1(\tilde X,\rho^*F_i)$. Proposition 3.16
implies that if $ker \rho^*_i\neq 0$ then $F_i=\Cal O(\tau)$.
Hence the kernel of the tangent map $\rho^*$ is trivial since it
follows from  lemma 3.27 that $H^1(X,F_i)=0$. This implies  that
$\rho^*$ is a local imbedding.\par

Part b) is a consequence of c), hence we first consider  c).
Theorem 3.26 states $\rho^*E$ has the decomposition $\rho^*E\cong
\tilde E_1\otimes\Cal O^k\oplus ...\oplus \tilde E_m\otimes \Cal
O^k$ with simple vector bundles $\tilde E_i$. If $\rho^* F \cong
\rho^*E$ then $\rho^*F$ inherits also a decomposition
$\rho^*F\cong \tilde F_1\otimes \Cal O^k\oplus ...\oplus \tilde
F_m\otimes \Cal O^k$ with $\tilde F_i\otimes \Cal O^k\cong \tilde
E_i\otimes \Cal O^k$. Since the $\tilde E_i$ are simple vector
bundles on $\tilde X$ it follows that $\tilde F_i\cong \tilde
E_i$. Also by theorem 3.26, we have a finite covering $p:X' \to X$
where $p^*E$ decomposes  as described in the theorem and equally
$p^*F$ decomposes into $p^*F\cong F'_1\otimes \Cal
O(\tau'_1)\oplus ...\oplus F'_m\otimes \Cal O(\tau'_m)$ with
${\rho'}^*F'_i=\tilde F_i$ It follows from lemma 3.24 that
$p^*E'_i\otimes \Cal O(\chi)\cong p^*F'_i$ for some character
$\chi:\pi_1(X')\to \Bbb C^*$, since ${\rho'}^*E'_i=\tilde
E_i\cong\tilde F_i= {\rho'}^*F'_i$ and $\tilde E_i$ is simple.
Hence c) follows from the decomposition for $p^*E$.
\par

To prove b) we first claim that there is a finite unramified
Galois covering $\hat p:\hat X \to X$ associated with $E$ such
that $\rho^*F\cong \rho^*E$ if and only if $\hat p^*F\cong \hat
p^*E$. The previous paragraph states that if $\rho^*F\cong
\rho^*E$ then:
$$p^*F\cong E'_1\otimes \Cal O(\tau_1)\oplus ... \oplus E'_m\otimes
\Cal O(\tau_m) \tag {3.3}$$ \noindent  where $\Cal O(\tau_i)$ are
flat bundles of rank $k$ associated with finite representations
(not the same as in the theorem 3.26). The variety of
representations of $\pi_1(X') \to GL(k,\Bbb C)$ for a fixed k,
$M(\pi_1(X'),k)$, is a finite set of points, since
 by lemma 3.27 $H^1(X',End\Cal O(\tau))=0$ for all  representations
$\tau$ and hence $M(\pi_1(X'),k)$ is zero dimensional. The
finiteness of the set of representations implies the existence of
a finite Galois cover $\hat p: \hat X \to X$ where $\hat p^*F\cong
\hat p^*E$ if $\rho^*F\cong \rho^*E$. The result follows then by
the lemma:

\proclaim {Lemma 3.28} Let $g:Y \to X$ be a finite unramified
Galois covering of $X$ and $E$ an absolutely stable bundle on $X$.
If $F$ is a vector bundle on $X$ such that $g*F\cong g^*E$ then
$F$ belongs to a finite collection of isomorphism classes of
vector bundles on $X$.
\endproclaim
\demo {Proof} If $g^*E$ is a simple vector bundle then the proof
of lemma 3.24 gives the result. More precisely, it shows that
$F\cong E\otimes \Cal O(\chi)$ where $\chi: G\to \Bbb C^*$ is a
character of the Galois group $G$ of the cover $f$.\par If $f^*E$
is not  simple applying the argument in theorem 3.26
 we get that $g^*E\cong E_1\otimes \Cal O^k\oplus ...\oplus E_m\Cal O_k$
 and $H^0(Y,Endg^*E)=\bigoplus^m_{i=1} M(k)$ for some $k$ dividing $\text
{rk} E$. The vector bundles $E$ and $F$ are quotients of two
different actions of the Galois group $G$ on $p^*G$. The quotient
of action of $G$ on $p^*E$ is up to isomorphism determined by the
isomorphism class of induced representation $\tau: G \to
GL(mk,\Bbb C)$. Our result follows  since the number of
isomorphism classes of  representations $\tau: G \to GL(mk,\Bbb
C)$ is finite.
 \hfill \ \qed
\enddemo

  \hfill \ \qed
\enddemo

\noindent Remark: We have a similar result for $H$-stable bundles
if $\tilde K_{eff}$ satisfies $P1$ or $P1'$.

\

What about the map of the space of all bundles (omitting the
discussion of wether it can be well defined)? Notice that for any
given filtration of saturated subsheaves in a vector bundle $V$
there is a blow up $X'$ of $X$ such that the pullback of this
filtration becomes a filtration of vector bundles (see Moishezon
[Mo69] lemma 3.5). In particular, for any vector bundle $V$ on $X$
one can use the Harder-Narasimhan filtration. Since the algebra of
holomorphic functions on $\tilde X$ does not change after changing
blowning up, any conclusion about the function theory for $\tilde
X'$ holds for $\tilde X$. It follows from the above that if $P1'$
holds then the pullback map for all bundles is non-injective
modulo representations of $\pi_1(X)$ only if there are cocycles
$\alpha \in H^1(X,V)$ such that $\rho^*\alpha=0$.

\

The following are some remarks about how to use the above results
to show that the universal cover of a projective variety has a
nonconstant holomorphic function.

\proclaim {Proposition 3.29} Let $X$ be a projective manifold of
dimension $n$ and $X'$ be an infinite unramified cover of $X$ then
$H^n(X', \Cal F)=0$ for any coherent sheaf $\Cal F$ on $X'$.
 \endproclaim

\demo {Proof} The result follows from Cech cohomology and Leray
coverings  if any noncompact cover of a n-dimensional projective
variety is covered by $n$ Stein open subsets. Pick $n-1$ generic
hyperplane sections $H_i$ and  let $C=H_1\cap ... \cap H_{n-1}$.
By Lefschetz theorem $C$ is a smooth curve such that $\pi_1(C)\to
\pi_1(X)$ is a surjection. This implies that the pre-image of $C$
in $X'$ is an irreducible  noncompact curve $C'$. Hence $C'$ is
Stein (Behnke-Stein theorem).  The infinite cover $X'$ is covered
by the pre-images $U_i$ of $X \setminus H_i$ in $X'$ and a
neighborhood of $C'$. The pre-images $U_i$ are Stein open subsets
of $X'$ since any unramified cover of a Stein manifold is Stein.
To conclude, $C'$ has  an open  Stein neighborhood in $X'$ since
$C'$ is a Stein closed subvariety of $X'$ (Siu [Si76]). \hfill \
\qed

\enddemo

\noindent Remark:  Proposition 4.26 implies that for surfaces the
structure of the space of the moduli space of vector bundles on
$\widetilde{X}$ should be similar to the structure of the moduli
space of vector bundles on a curve. Namely the groups $H^2(\tilde
X,\Cal F)$ vanish for any coherent sheaf $\Cal F$. In particular,
there are no algebraic obstructions in $H^2(\tilde X,End E)$ to
deform a vector bundle $E$ along a cocycle in $H^1(\tilde X,EndE)$
though there may be an analytic one (problem of convergency). We
expect  that any bundle of rank $\geq 2$ has a complete flag of
subbundles if there is a complete flag of topological subbundles.
This would imply that the K-group $K_{0}(\widetilde{X})$ reduces
to $\operatorname{Pic}(X)\times {\Bbb Z}$. The above motivates the
authors' expectation  that many different bundles on $X$ coincide
after pulling back to $\widetilde{X}$.

\

\

\head 4. Geometric vanishing theorem for negative bundles
\endhead

\

\

The arguments used in the proof of Theorem  A can be used to give
an alternative proof of the vanishing theorem for negative vector
bundle $V$ over a projective manifold $X$ whose $\operatorname{rk}
V < \dim X$.

\proclaim {Theorem 4.1} If $V$ is a negative vector bundle on a
projective manifold $X$ with $\operatorname{rk} V < \dim X$, then
$H^1(X,V) = 0$.
\endproclaim

\demo {Proof:} Suppose it exists a nontrivial $s \in H^1(X,V)$ and
let:

$$0\to V\to V_s \to\Cal O\to 0$$

\noindent be the associated extension.  As in the Theorem A,
consider the dual exact sequence and $A_s= \Bbb P(V_s^*) \setminus
\Bbb P(V^*)$ be an affine bundle, which by the negativity of $V$
is strictly pseudoconvex. Let $r:A_s \to St(A_s)$ be the Remmert
reduction, where $r$ is proper contracting $M=\cup_{i=1}^k M_i$
and $St(A_s)$ is a Stein space with isolated singularities.\par

The aim is to obtain a contradiction from topological conditions.
The Stein space $St(A_s)$ has $\dim _{\Bbb C} St(A_s)=\dim_{\Bbb
C}X +r$ and hence it has the homotopy type of a simplicial complex
of real dimension at most equal to $\dim_{\Bbb C}X +r$. On the
other hand, $St(A_s)=A_s/(\coprod_iM_i)$ as a topological space
and so for the reduced singular homology of $A_s$
$\widetilde{H}_{i}(St(A_s), {\Bbb C}) = H_{i}((A_s,\coprod_{i}
M_{i}), {\Bbb C})$. Now the long exact homology sequence of the
pair $(A_s,\coprod_{i} M_{i})$ together with the fact that
$\coprod_{i} M_{i}$ is compact of complex dimension strictly less
than $\dim_{\Bbb C}X=n$ (by proposition 2.5) gives that
$H_{2n}(A_s,{\Bbb C}) \cong \widetilde{H}_{2n}(St(A_s),  {\Bbb C})
= H_{2n}(St(A_s), {\Bbb C})$.\par

In conclusion, $St(A_s)$ as a Stein manifold of $\dim_{\Bbb
C}St(A_s)=n+r<2n$ must have $H_{2n}(St(A_s),\Bbb C)=0$. The
previous argument gives $H_{2n}(St(A_s), {\Bbb C}) \cong
H_{2n}(A_s,{\Bbb C})$. The contradiction follows since $A_s$ as an
affine bundle over $X$  is homotopicaly equivalent to $X$ and
therefore $H_{2n}(A_s,{\Bbb C}) \cong H_{2n}(X,{\Bbb C})\neq 0$.
\enddemo

\

\noindent Remark: This proof works also for normal projective
varieties.

\

\

\end